\documentclass[a4paper, 11pt]{article}

\usepackage{amssymb, amsmath}
\usepackage{epsfig}
\usepackage{float}
\usepackage{subfig}
\usepackage[margin=3cm]{geometry}
\usepackage{graphicx}
\usepackage{times}
\usepackage{url}
\usepackage{xcolor}
\usepackage{xspace}

\colorlet{dblue}{blue!80!black}
\colorlet{dgreen}{green!50!black}

\colorlet{darkorange}{orange!70!black}

\newcommand{\pforest}{\texttt{p4est}\xspace}
\newcommand{\hypre}{\texttt{hypre}\xspace}

\newcommand{\eqnlab}[1]{\label{eqn:#1}}
\newcommand{\eqnref}[1]{\eqref{eqn:#1}}

\newcommand{\seclab}[1]{\label{sec:#1}}
\newcommand{\secref}[1]{Section~\ref{sec:#1}}

\newcommand{\tablab}[1]{\label{tab:#1}}
\newcommand{\tabref}[1]{Table~\ref{tab:#1}}

\newcommand{\figlab}[1]{\label{fig:#1}}
\newcommand{\figref}[1]{Figure~\ref{fig:#1}}

\newcommand{\eg}{e.g.\ }
\newcommand{\ie}{i.e.,\xspace}

\newcommand{\field}[1]{\mathbb{#1}}

\newcommand{\R}{\field{R}}

\newcommand{\hdiv}{\ensuremath{\boldsymbol{H}(\dive, \Omega)}}
\newcommand{\ltwo}{\ensuremath{L}^{2}(\Omega)}
\newcommand{\sob}[1]{\ensuremath{H}^{#1}(\Omega)}
\newcommand{\dirich}{\ensuremath{H_{0,D}^{1}(\Omega)}}
\newcommand{\neum}{\ensuremath{\boldsymbol{H}_{0,N}(\Omega)}}
\newcommand{\rtzero}{\ensuremath{\mathcal{RT}_0}\xspace}

\newcommand{\K}{\ensuremath{\boldsymbol{\mathcal{K}}}}
\newcommand{\Kinv}{\ensuremath{\boldsymbol{\mathcal{K}^{-1}}}}
\newcommand{\diag}[1]{\ensuremath{\mathrm{diag}(#1)}}

\DeclareMathOperator{\idd}{\mathrm{Id}}
\DeclareMathOperator{\grad}{\mathrm{grad}}
\DeclareMathOperator{\dive}{\mathrm{div}}
\newcommand{\diver}[1]{\ensuremath{\dive\boldsymbol{#1}}}

\newcommand{\ev}{\,\cdot\,}
\newcommand{\norm}[2]{\lVert#1\rVert_{#2}}
\newcommand{\seminorm}[2]{\lvert#1\rvert_{#2}}
\newcommand{\dual}[2]{\ensuremath{\langle#1,#2\rangle}}

\newcommand{\bm}{\boldsymbol}

\newcommand{\ds}{\ensuremath{\,\mathrm{d}s}}

\newcommand{\dx}{\ensuremath{\,\mathrm{d}x}}

\newcommand{\trhcm}{
Authors B.\ and F.\ gratefully acknowledge financial
support by the SFB/TR 32 ``Patterns in Soil-Vegetation-Atmosphere Systems:
Monitoring, Modeling, and Data Assimilation'' funded by the Deutsche
Forschungsgemeinschaft (DFG), as well as travel support from the DFG-funded
Hausdorff Center for Mathematics, Bonn, Germany.%
}

\newcommand{\sfb}{
Author M.\ acknowledges the support of the former SFB611
``Singular Phenomena and Scaling in Mathematical Models'' funded by the
Deutsche Forschungsgemeinschaft (DFG).%
}

\title{An AMG saddle point preconditioner with application to
mixed Poisson problems on adaptive quad/cube meshes}
\author{Carsten Burstedde, Jose A.\ Fonseca, and Bram Metsch}
\begin{document}

\maketitle

\begin{abstract}
We investigate various block preconditioners for a low-order Raviart-Thomas
discretization of the mixed Poisson problem on adaptive quadrilateral meshes.
% We compare the effectiveness of various block preconditioners for saddle point
% problems arising from a mixed element (MFE) discretization of the Poisson
% equation, incorporating a new block-AMG saddle point approach.
% We consider lowest order Raviart-Thomas elements on adaptive
% quadrilateral meshes.
% The first preconditioner is a lumped diagonal
% approximation of the system matrix.  The second employs an exact block
% factorization of the system matrix combined with a lumped diagonal
% approximation of the velocity term and an algebraic multigrid (AMG) V-cycle on
% an approximate Schur complement.
% The third one is a dedicated
% AMG solver for saddle point problems (SPAMG).
In addition to standard diagonal and Schur complement preconditioners, we
present a dedicated AMG solver for saddle point problems (SPAMG).
A key element is a stabilized prolongation operator that couples the flux and
scalar components.
Our numerical experiments in 2D and 3D show that the SPAMG preconditioner
displays nearly mesh-independent iteration counts for adaptive meshes and
heterogeneous coefficients.
%, a property that is lost with the common Schur complement approach.
\end{abstract}

%\jose{I removed the [H] specifiers from figures.
%      Please make sure that all figures are referenced correctly and that their
%      captions contain everything the reader needs to understand them, with
%      references from the caption into the text to point to equation and
%     section numbers.}

\section{Introduction}
\seclab{intro}

In many applications one is more interested in the gradient of the solution
of the Poisson equation than in the solution itself.
One practical example is solving the Richards equation for subsurface flow,
where the scalar variable represents the water pressure \cite{Richards31, Bear13}.
One approach is
to use standard methods such as finite differences (FD), finite elements (FE) or
finite volumes (FV), to obtain an approximate solution and then to compute itself
gradient via numerical differentiation, which may lead to a loss in accuracy
\cite{Arnold90}.
A mixed finite element (MFE) discretization then appears to be a more natural
choice since the gradient of the solution is part of the unknown variables of
the formulation.
Generally, in MFE methods the vector valued quantity is approximated at least
with the order of accuracy of the scalar unknown \cite{BrezziFortin91}.
In addition, the computed solution satisfies mass
conservation at the element level, a property highly relevant for simulating
fluid flow problems where the governing partial differential equations (PDE)
are derived from mass balance laws \cite{BoffiBrezziFortin13}.

MFE discretizations lead to symmetric indefinite systems of algebraic equations
that may be solved by iterative methods. In order to obtain a solution with
a reasonable investment of computational resources, the use of
optimal preconditioners becomes mandatory. Although the existence
of an optimal preconditioner for general (anisotropic) coefficients remains an
open question, there has been a significant amount of work in this direction.
\cite{ArnoldFalkWinther97} proposes two preconditioners for the
operator $(I - \grad \dive)$ in two dimensions, one based on domain
decomposition and another on multigrid.
The latter has been generalized to operators of the form  $(\rho I - \mu \grad
\dive)$ and three-dimensional problems \cite{ArnoldFalkWinther00}.
From the observation that a MFE discretization of a generali\-zed diffusion equation
is well posed in a product of two different discrete function spaces,
\cite{PowellSilvester03} and
\cite{Powell05} propose two block-diagonal preconditioners. One of them is what we
will refer to as Schur complement preconditioner and consists
of a lumped diagonal approximation of the $(1,1)$ block and an algebraic
multigrid V-cycle to approximate a Schur complement on the $(2,2)$ block.
The preconditioners introduced in these papers are shown to be optimal with
respect to the mesh size and various classes of coefficients (such as a
conductivity tensor).
Recent work \cite{KrausLazarovLymberyEtAl16} introduces a new preconditioner
whose key ingredient is the approximation of the $(1,1)$ block with an
auxiliary space multigrid method \cite{KrausLymberyMargenov15} that offers
optimality with respect to a larger class of coefficients.
The authors include numerical evidence for two dimensional problems posed on
uniform rectangular meshes.

For some problems we may require a very fine mesh in order to correctly resolve
the phenomena we are trying to model.
A uniform mesh might be undesirable or even impractical given the computational
resources it requires.
One solution to this problem is to use locally refined meshes, which use the
correct resolution only in the portion of the domain that really requires it
\cite{BabuskaRheinboldt78, BergerColella89, Dorfler96}.
There are essentially two possibilities:
either we allow the so called hanging nodes or not.
In the context of MFE, implementations making use of adaptively refined meshes
with hanging nodes are not the standard, nevertheless the theoretical framework
has been established in the early 90s \cite{EwingLazarovRussellEtAl90,
EwingWang92}.
Unfortunately, none of the above mentioned methods for preconditioning address
the case of adaptively refined meshes.
Hence, the goal of this paper is to present a multigrid preconditioner that
retains its robustness in such a case for both two and three dimensional problems.
In addition, we aim to allow for variable matrix-valued coefficients.
In our approach, we acknowledge the different scales in the saddle point
structure and introduce a prolongation operator that suitably couples the
vector and scalar unknowns.
We demonstrate superiority of this approach to the  Schur complement method
using a variety of numerical examples.

% \todo{Brief run through the ideas that go into this paper.}
% In section 2 we review the strong and weak formulation of the problem, following
% by the discretization and the local refinement approach chosen for this document.
% In section 3 we introduce an algebraic multigrid method suitable for saddle
% point problems and discuss two different kind of smoothers for such problems. Finally,
% we evaluate the effectiveness of this multigrid method as a preconditioner compared
% to the classical Schur complement in the framework  of uniform and adaptive meshes.

\section{Problem formulation}
\seclab{mixed}

In this section we briefly review the % strong and weak
mixed formulation
of a Poisson-type partial differential equation.
The material covered here is standard and found in many textbooks; see \eg
\cite{BoffiBrezziFortin13, BrennerScott02}, augmented with a few recent
results.

Let $\Omega\in\R^{d}$ for $d\in\{2,3\}$ be a bounded Lipschitz domain
with boundary $\Gamma=\partial\Omega$.
$\sob{m}$ will denote the standard Hilbert space of functions in $\ltwo$
whose weak derivatives up to order $m\geq0$ are also square integrable.
$\sob{m}$ is endowed with the usual norm and seminorm,
\begin{equation}
\eqnlab{hnorm}
\norm{v}{m}^2 := \sum_{|\alpha|\leq m}\int_{\Omega}|D^{\alpha}v|^2 \dx ,
\quad
\seminorm{v}{m}^2 := \sum_{|\alpha|= m}\int_{\Omega}|D^{\alpha}v|^2 \dx
.
\end{equation}
We define the velocity space
\begin{equation}
\eqnlab{hdiv}
\hdiv:=
\{\bm{v}\in(\ltwo)^{d} : \diver{\bm{v}}\in\ltwo \}
,
\end{equation}
which is a Hilbert space with the norm
\begin{equation}
\norm{\bm{v}}{\dive} ^2
=
\norm{\bm{v}}{0}^{2} +
\norm{\diver{\bm{v}}}{0}^{2}
.
\end{equation}
We assume that
$\Gamma=\Gamma_{D}\cup\Gamma_{N}$,
with $\Gamma_{D}\cap\Gamma_{N}=\emptyset$,
and that $\Gamma_{D}$ has nonzero $(d-1)$-dimensional Lebesgue
measure. Lastly, define (in the weak sense)
\begin{equation}
\eqnlab{hdirichlet}
\dirich :=
\{\phi\in\sob{1} : \phi|_{\Gamma_{D}}=0 \}
.
\end{equation}
In the following, we use bold mathematical symbols to denote
vectors and matrices over $\R$ and the normal font for scalar quantities.

We consider the equation
\begin{subequations}
\eqnlab{poisson:strong}
\begin{align}
  -\diver{\K}(\bm x) \nabla p = f &\quad\text{ in } \Omega, \\
  p = p_0 &\quad\text{ on }\Gamma_D,\\
  \K(\bm x)\nabla p \cdot \bm n  = g &\quad\text{ on } \Gamma_{N}
  , \eqnlab{poisson:strong:neumann}
\end{align}
\end{subequations}
where $\bm n$  is the outer normal vector to the boundary $\Gamma$.
The conductivity tensor $\K(\bm x)$ is
a $d\times d$ symmetric positive definite matrix whose smallest eigenvalue
is bounded uniformly away from zero. Furthermore, the data is required to satisfy
\begin{subequations}
\eqnlab{poisson:data}
\begin{align}
  f\in\ltwo, \quad p_{0} \in H^{1/2}(\Gamma_{D})
  \quad\text{and}\quad g\in L^{2}(\Gamma_{N})
  .
\end{align}
\end{subequations}
Introducing the variable $\bm u=\K \nabla p$ leads to the  mixed first-order
system
\begin{subequations}
\eqnlab{poisson:system}
\begin{align}
  \bm u  = \K \nabla p  &\quad\text{ in } \Omega, \\
  -\diver{\bm u}  = f  &\quad\text{ in } \Omega, \\
  p = p_0 &\quad\text{ on } \Gamma_{D},\\
  \bm u \cdot \bm n = g &\quad\text{ on } \Gamma_{N}
  .
  \eqnlab{poisson:system:neumann}
\end{align}
\end{subequations}

\subsection{Weak formulation}

To derive the mixed weak formulation of \eqnref{poisson:system} we
introduce the following space,
\begin{equation}
  \neum
  :=
  \{\bm\tau\in\hdiv : \dual{\bm{\tau}\cdot\bm{n}}{\phi}= 0
    \text{ for all } \phi\in\dirich\}
  .
\end{equation}
The dual paring $\dual{\ev}{\ev}$ is defined via Green's formula;
see \cite[pp.~50]{BoffiBrezziFortin13}.
Multiplying the first equation in \eqnref{poisson:system} by $\Kinv$ and
then by a test function $\bm{v} \bm\in \neum$, the second
by some $q\in\ltwo$, integrating over $\Omega$ and
using Green's formula on the gradient term yields
the following weak formulation: Find $(\bm u,p)\in  \neum\times\ltwo$
such that
\begin{subequations}
\eqnlab{poisson:integral}
\begin{align}
  \int_{\Omega} \Kinv \bm u\cdot\bm{v} \dx +
  \int_{\Omega} p\,\diver{\bm{v}} \dx &=
  \int_{\Gamma_{D}} p_0(\bm{v}\cdot\bm n) \ds
  \quad\text{for all } \bm{v}\in \neum, \\
  \int_{\Omega}q\,\diver{u}\dx &=
  -\int_{\Omega}fq\dx
  \quad\text{for all } q\in\ltwo
  .
\end{align}
\end{subequations}
Let us define the bilinear forms
$a:\hdiv\times\hdiv\to\R$ and $b:\hdiv\times\ltwo\to\R$,
\begin{equation}
\eqnlab{bforms}
  a(\bm u,\bm v):= \int_{\Omega} \Kinv \bm u\cdot\bm{v} \dx, \quad
  b(\bm v,p):=\int_{\Omega} p\,\diver{\bm{v}}
  ,
\end{equation}
and linear functionals $\mathbb{A}:\hdiv\to\R$,  $\mathbb{B}:\ltwo\to\R$,
\begin{equation}
\eqnlab{lforms}
   \mathbb{A}(\bm{v}):= \int_{\Gamma_{D}} p_0(\bm{v}\cdot\bm n) \ds, \quad
   \mathbb{B}(q):=-\int_{\Omega}fq\dx
  .
\end{equation}
To enforce the Neumann boundary condition \eqnref{poisson:strong:neumann},
\eqnref{poisson:system:neumann}, let $\tilde{p}$ a classical solution
of \eqnref{poisson:strong} with $f=0$ and $p_0=0$.
Then, taking $\bm{\tilde{u}} = \K\nabla\tilde{p}$ leads to the following mixed
weak formulation of \eqnref{poisson:system}:
Find $\bm{u}=\bm{\tilde{u}}+\bm{u}_{0}$ with
$\bm{u}_{0}\in \neum$ and $p\in\ltwo$ such that
\begin{subequations}
\eqnlab{mixed}
\begin{alignat}{2}
  a(\bm{u}_0,\bm{v}) + b(\bm{v},p) &=
  \mathbb{A}(\bm{v})-a(\bm{\tilde{u}},\bm{v})
  := \tilde{\mathbb{A}}(\bm v)
  \qquad&\text{for all } &\bm{v}\in \neum
  , \\
  b(\bm{u},q) &= \mathbb{B}(q)-b(\bm{\tilde{u}},q)
  := \tilde{\mathbb{B}}(q)
  &\text{for all } &q\in\ltwo
  .
  \eqnlab{neumann_enforce}
\end{alignat}
\end{subequations}
Existence and uniqueness of a solution for the
problem \eqnref{mixed} follows from standard arguments, namely
establishing coercivity  of the bilinear form $a(\ev,\ev)$, and the
Lady\v{z}enskaja-Babu\v{s}ka-Brezzi (LBB) condition;
see \eg \cite{BrezziFortin91}.

\subsection{Discretization}
\seclab{discret}

We pick finite dimensional subspaces $\bm{X}_{0}^{h}\subset \neum$
and $V^{h}\subset\ltwo$ and define the following problem:
Find $(\bm u_{h},p_{h})\in \bm{X}_{0}^{h}\times V^{h}$ such that
\begin{subequations}
\eqnlab{mixed:discrete}
\begin{alignat}{2}
  a(\bm{u}_{h},\bm{v}_{h}) + b(\bm{v}_{h},p_{h}) &= \tilde{\mathbb{A}}(\bm{v}_{h})
  \quad&\text{for all } &\bm{v}_{h}\in \bm{X}_{0}^{h}
  , \\
  b(\bm{u}_{u}, q_{h}) &= \tilde{\mathbb{B}}(q_{h})
  \quad&\text{for all } &q_{h}\in V^{h}
  .
\end{alignat}
\end{subequations}
% %\jose{Please make sure that the coefficients in $b()$ are always in the same order!
% %      This is inconsistent here.}
% where the linear forms $\tilde{\mathbb{A}}$ and $\tilde{\mathbb{B}}$ take into
% account the Neumman boundary condition like in \eqnref{mixed}, that is
% \begin{subequations}
% 	\begin{align}
% 	 \tilde{\mathbb{A}}(\bm v) &:=
% 	\begin{cases}
%         \mathbb{A}(\bm v)-a(\bm{\tilde{u}},\bm{v}) &\text{ if }
% 		\Gamma_N\neq\emptyset , \\
% 		\mathbb{A}(\bm v) &\text{ else,}
% 	\end{cases} \\
% 	\tilde{\mathbb{B}}(q)&:=
% 	\begin{cases}
% 		\mathbb{B}(q)-b(\bm{\tilde{u}},q) &\text{ if } \Gamma_N\neq\emptyset,\\
% 		\mathbb{B}(q) &\text{ else.}
% 	\end{cases}
% 	\end{align}
% \end{subequations}
To ensure that \eqnref{mixed:discrete} is well posed, the pair
of spaces $(\bm{X}_{0}^{h}, V^{h})$
%\todo{What is $X_g^h$? Typo, should be X_0^h}
must be chosen such that the LBB condition is fulfilled for the discrete
problem. Many spaces with this property have been developed since the early
seventies, such as the Raviart-Thomas  \cite{RaviartThomas77} and
Brezzi-Douglas-Marini \cite{BrezziDouglasMarini85} spaces.
In this paper we choose the lowest order Raviart-Thomas discretization
\rtzero defined on rectangles/hexahedra.
Hence, in an element $\Omega_e$ the velocity and pressure test functions take
the form
\begin{equation}
  \eqnlab{rtzero}
  \bm v |_{\Omega_e} =
  \begin{Bmatrix}
    (a_0 + b_0x, a_1 + b_1y )^{T} &\quad\text{if } d=2\\
    (a_0 + b_0x, a_1 + b_1y, a_2 + b_2z )^{T} &\quad\text{if } d=3
  \end{Bmatrix}
  ,\quad
  p |_{\Omega_e} = c_0
  ,
\end{equation}
respectively, where $a_i, b_i$ for $i=0,\ldots,d - 1$ and $c_0$ are constants.
The degrees of freedom are shown in \figref{rtdof}.
\begin{figure}
  \centering
  \subfloat[]{\includegraphics[scale=1.2]{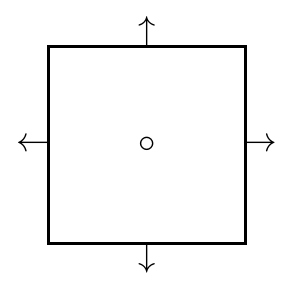}}
  \hspace{4em}
  \subfloat[]{\includegraphics[scale=2]{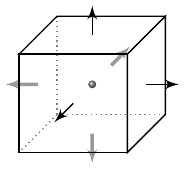}}
  \caption{Dregrees of freedom for the rectangular \rtzero element in two (a) and
	three (b) dimensions. For the velocity, the degrees of freedom are normal
	components at the edge (face) mid sides of an element.
	The pressure is located at the center of the element.}
  \figlab{rtdof}
\end{figure}
The following approximation properties are well known for \rtzero in the
context of affine elements defined on uniform meshes \cite{BoffiBrezziFortin13},
\begin{subequations}
\eqnlab{conv_rate}
\begin{align}
  \norm{\bm{u}-\bm{u}_{h}}{0} &\leq Ch\norm{\bm{u}}{1},\\
  \norm{p-p_{h}}{0} &\leq c h\left(\norm{p}{1} + \norm{\bm u}{1}\right)
  ,
\end{align}
\end{subequations}
where we assume that the pair $(\bm u, p)$ fulfils the regularity requirements
required by the right hand side of \eqnref{conv_rate}.

\subsection{Adaptivity}

For well behaved (smooth) problems posed on convex domains,
the use of a uniform mesh usually offers
satisfactory results  when computing a numerical solution, that is,
there is an optimal trade-off between numerical effort (computational resources)
invested and effective reduction of the error.
Nevertheless,
there are situations in which the mesh resolution required
to accurately reproduce the physical behavior of the underlying PDE
becomes computationally impractical if imposed on the whole domain.

Adaptive mesh refinement (AMR) provides a valuable tool in order to
reduce the computational complexity in such situations by increasing
the mesh resolution only locally where is required. As stated in section
\secref{discret}, we will work with meshes composed of rectangles/hexahedra.
The refinement schemes used in this paper include the case of a mesh with
hanging nodes, that is, we allow (a nonempty) intersection of two
elements to be a complete side of a neighboring element. Additionally, we
choose to use 2:1 balanced meshes: The length ratio between
a coarse and a fine element is at most of factor two; see \figref{amr}.
Non-balanced meshes would be possible as well but require more technical work
regarding the definition of MFE spaces and parallelization.
\begin{figure}
  \centering
  \includegraphics[scale=1.3]{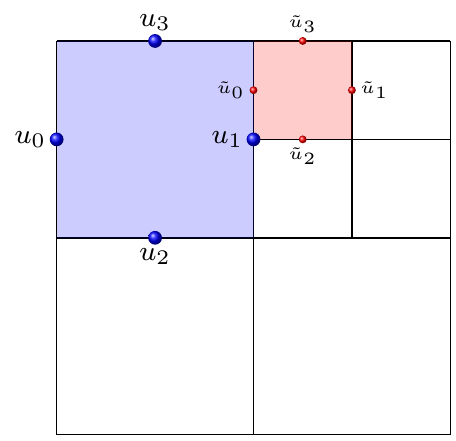}
  \caption{Locally refined mesh with hanging nodes. To enforce continuity
	of the flux normals for a \rtzero discretization, the
	velocity value on a hanging node is defined by the
	corresponding non-hanging node lying in the same edge.
	For the case in the picture we  define $\tilde{u} _0 := - u_1$.
%        \jose{Edited---is this our sign convention, or are we
%              choosing always positive (flux from left to right)?
%              We are choosing it always positive}
          }
  \figlab{amr}
\end{figure}
Continuity of fluxes across interfaces for this kind of meshes can be enforced
in several ways.
One is to eliminate te degrees of freedom at hanging nodes
\cite{EwingLazarovRussellEtAl90}.
Another approach is to use Mortar finite elements
\cite{ArbogastCowsarWheelerEtAl00}.
We will
follow the former.  Due to our assumption of 2:1 balance, for a
\rtzero discretization given a hanging node with flux value $\tilde{u}$, there
is only one non-hanging node $u$ lying in the same edge/face; see \figref{amr}.

Estimates for locally refined meshes using hanging nodes have
been studied in \cite{EwingLazarovRussellEtAl90, EwingWang92}.
Essentially, it is shown that the \rtzero spaces still respect the LBB
condition after introducing locally refined grids, and hence the estimates
\eqnref{conv_rate} remain valid.

\subsection{Preconditioning}

The discretization of \eqnref{poisson:system} via stable MFE
leads to a saddle point problem defined by the following block matrix,
\begin{equation}
  \eqnlab{block}
  \mathcal{A}_h:=
  \begin{bmatrix}
  \bm{A}  & \bm{B}^{\top} \\
  \bm{B} & 0 \\
  \end{bmatrix}
  \sim
  \begin{bmatrix}
  \idd_h  & - \grad_h \\
  \dive_h &  0 \\
  \end{bmatrix}
  ,
\end{equation}
where $\bm{A}$ is the vector mass matrix and $\bm{B}$ is the
discrete divergence operator. This is a well studied
problem and there are several solution methods available; see
\eg \cite{BenziGolubLiesen05} for a comprehensive review.
They range from Uzawa algorithms and its variants
\cite{BramblePasciakVassilev97}, projection methods \cite{GuermondMinevShen06},
to block factorization methods \cite{ElmanSilvesterWathen14, MardalSundnesLangtangenEtAl03}.
In contrast to most methods that treat flux and pressure individually, we
introduce a monolithic multigrid method originally developed for the Stokes
equations \cite{Metsch13} that is to our knowledge yet unpublished.

To prepare the following exposition, let us briefly discuss notation and some
background. We consider the case that the matrix \eqnref{block} is
symmetric and indefinite.
The factorization
\begin{equation}
  \eqnlab{block:fact}
  \mathcal{A}_h=
  \begin{bmatrix}
  \bm{I} &  0 \\
  \bm{BA^{-1}} & \bm{I} \\
  \end{bmatrix}
  \begin{bmatrix}
  \bm A & 0 \\
   0 & - \bm{S} \\
  \end{bmatrix}
  \begin{bmatrix}
  \bm{I}  &
  \bm{A^{-1}B^{\top}}\\
   0 &
  \bm{I} \\
  \end{bmatrix},
  \qquad
  \bm{S} = \bm{BA^{-1}B^{\top}}
  ,
\end{equation}
implies that $\mathcal{A}_h$ is congruent to a block diagonal matrix \cite{Fiedler08}.
This fact motivates the use of a preconditioner of the form
\begin{equation}
  \eqnlab{precond}
  \mathcal{B}=
  \begin{bmatrix}
  \bm{M} & 0 \\
  0 & - \bm{N} \\
  \end{bmatrix}
  ,
\end{equation}
where $\bm M$ and $\bm N$ satisfy
\begin{subequations}
  \begin{align}
    \bm{M}\bm{A} \approx \bm{I}, \eqnlab{precond:properties:a}\\
    \bm{N}\bm{S} \approx \bm{I} \eqnlab{precond:properties:b}
    .
  \end{align}
\end{subequations}
A simple choice is to take $\bm{M}$ as the inverse of the lumped mass matrix $\bm{A}$.
The Schur complement $\bm{S}$ represents the operator $-\Delta_{h}$. Hence,
\eqnref{precond:properties:b} suggests $\bm{N}\approx \Delta_h^{-1}$, and the second
block of the preconditioner should approximate the inverse of a discrete
Laplacian in pressure space.
Then, we can use a solver for elliptic operators such as multigrid to
apply $\bm{N}$. Nevertheless, it is a concern that our pressure belongs to $\ltwo$, and
strictly speaking we do not have the required regularity to apply $\Delta^{-1}_h$.
An option to deal with this problem is to use the auxiliary space technique
\cite{Xu96}, in which the idea is to use a multigrid preconditioner
for continuous pressure elements and then project it in to the desired
space of discontinuous pressure elements in combination with a suitable
smoothing operator. We note that some approaches only use a one-sided factorization of
\eqnref{block:fact}, which leads to a block-triangle form of $\mathcal B$;
\eg \cite{ElmanSilvester96,KronbichlerHeisterBangerth12}.
For problems in which the (1,1) block of $\mathcal{A}_h$ has a non-symmetric
term, this variant offers a faster convergence of the iterative solver at the
price of evaluating an additional sparse matrix vector product compared to
the block diagonal preconditioner; see \cite{FischerRamageSilvesterEtAl98}.

In the context of the (Navier-)Stokes equations, dedicated approximations
to the inverse of the Schur complement have been proposed
\cite{ElmanSilvester96, Elman99}.
These include the pressure Schur complement methods
\cite{Turek99} and the least squares commutator preconditioner a.k.a.\ BFBt
and its extensions \cite{ElmanSilvesterWathen14}.
%, RudiStadlerGhattas17}.
Following
the presentation from \cite{RudiStadlerGhattas17}, the BFBt approximation
of the inverse of the Schur complement can be written
\begin{equation}
  \eqnlab{bfbt}
  \bm S_{BFBt}^{-1} :=
	(\bm{BC}^{-1} \bm B^\top)^{-1}
	(\bm{BC}^{-1}\bm{AD}^{-1}\bm B^\top)
	(\bm{BD}^{-1}\bm B^\top)^{-1}
  ,
\end{equation}
where  $\bm C$ and $\bm D$ are diagonal and symmetric positive definite matrices.
The original choice sets $\bm C$ and $\bm D$ to the
lumped velocity mass matrix of the system.
New modifications have been
introduced in an effort to improve the effectiveness of the preconditioners in
cases where the equations present high variability in the (scalar) coefficients
\cite{RudiMalossiIsaacEtAl15, RudiStadlerGhattas17}.
For Poisson instead of Stokes, $\bm A$ is the velocity mass matrix instead of a
discrete Laplacian.
Hence, with $\bm C = \bm D \approx \bm A$ the method reduces to
$\bm S_{BFBt} \approx \bm S$,
% which reduces \eqnref{bfbt} to
the usual Schur complement.
% approximation.
% that
% We denote $\bm N := (\bm{B\tilde{A}}^{-1} \bm B^\top)^{-1}$,
% With $\bm{\tilde{A}} := \text{lump}(\bm A)$, thus
% .
% We approximate the action of $\bm N$
% with a multigrid method even though the \rtzero pressure is piecewise constant.

\section{Multigrid for Saddle Point Problems}
\seclab{SPAMG}
Multigrid (MG) methods provide efficient preconditioners for important
classes sparse of linear systems, in particular if the matrix system arises
from the discretization of an elliptic PDE. Their main advantage is that they scale
linearly in the number of unknowns $N$, \ie that they require only $O(N)$
computational work and memory.
Multigrid has been primarily developed for symmetric positive
M-matrices as they typically arise from FD/FV/FE discretizations of
(scalar) second order elliptic PDEs.

As already indicated in
the previous section, multigrid algorithms can be used inside Schur
complement preconditioners to invert one (or both, depending on the
application) of the diagonal blocks.
While this approach allows for an easy reuse of existing and
well-established techniques, it does no longer
guarantee linear convergence for the block system as a whole.  The
multigrid cycles are applied only to the sub-problem(s), while
couplings between the unknowns are handled by the outer iteration
(for example GMRES, BiCGStab or (inexact) Uzawa).  In consequence, the
outer iteration essentially determines the overall convergence speed,
even if the inner sub-blocks can be solved quickly.

The question is whether it is possible to build a multigrid hierarchy
for the coupled system to take into account the cross-couplings on all
levels. Several developments have been made in this
direction. Geometric multigrid methods for (Navier-)Stokes have been
proposed in \eg \cite{Vanka:1986,SchoeberlZulehner03}. In
\cite{Wabro:2003,Wabro:2004,Wabro:2006}, this has been extended to a
``semi-algebraic'' AMGe, where the coarse levels are still
determined geometrically. In this section, we will present the fully algebraic
method from the original thesis \cite{Metsch13}.
As in classical AMG for M-matrices, this method only requires the matrix for
building a robust multigrid hierarchy and hence can adapt itself to
difficulties such as anisotropic or jumping coefficients as well as non-uniform
meshes that cannot be coarsened geometrically.

\subsection{Algebraic multigrid (AMG)}
\label{section:amg_elliptic}
%\bram{Feel free to add to this paragraph.  Mind that this paper is on Poisson,
%not Stokes, and adapt discussion accordingly.  You may add your own bib file if
%you like}.
The whole hierarchy of grids $\left\{\Omega_l \right\}_{l=1}^L$, differential
operators $\left\{ {\bm \Lambda}_l \right\}_{l=1}^L $, and transfer operators
$\left\{{\bm P}_l \right\}_{l=1}^{L-1}$, $\left\{ {\bm R}_l \right\}_{l=1}^{L-1} $
needed in the multigrid cycle (the solution phase)
is computed from the system matrix ${\bm \Lambda}$ during a setup phase.
In this document we will understand the term ``grid''
as a set of indices, since AMG methods require no geometric mesh information.
The strength of AMG is its ability
to deal with difficulties in the operator, such as  heterogeneous,
strongly varying coefficients as well as unstructured meshes, for which a
hierarchy cannot be (easily) identified.
To set the stage for our proposed method, let us
briefly recapitulate the classical Ruge-St\"uben AMG setup
\cite{RugeStueben87, Stueben99}.
We start on the finest level $l=1$ using the
fine system matrix ${\bm \Lambda}_1 := {\bm \Lambda}$ and the finest set of degrees
of freedom $\Omega_1 = \{1,\ldots,n_1\}$.
%\footnote{} $$.
\begin{enumerate}
\item We decompose the grid $\Omega_l$ into the set of {\em fine grid
    points} $F_l$ and {\em coarse grid points} $C_l$. The latter form
  the next coarser grid $\Omega_{l+1} := C_l$ of size $n_{l+1}$.
\item  We compute the interpolation matrix
\begin{equation}
  \bm{P}_l \in \mathbb{R}^{n_l \times n_{l+1}}, \quad \bm{P}_l
  := \begin{bmatrix} \bm{P}_{l,F} \\ \bm I_{l,C} \end{bmatrix}.
\end{equation}
 The submatrix $\bm{P}_{l,F}$ contains the interpolation weights for
 all fine grid points $i \in F_l$. For the coarse grid
points $i \in C_l$, interpolation is just the identity.
\item We compute the next coarser matrix by the Galerkin product
  \begin{equation}
  \eqnlab{galerkinproduct}
    \bm{\Lambda}_{l+1} := \bm{P}_l^{\top} \bm{\Lambda}_l \bm{P}_l .
  \end{equation}
\end{enumerate}
The algorithm is then recursively applied to the input matrix
$\bm{\Lambda}_{l+1}$. If $n_l$ is small enough such that
$\bm{\Lambda}_l$ can be efficiently
factored by a direct solver, the recursion is terminated. In our
experiments we use a redundant serial solver if $n_l \le 1000$.

% We now discuss some aspects of the setup phase.
We recall
that for a fast multigrid algorithm, the error components not
efficiently reduced by smoothing (the smooth vectors) must be
well represented within the range of the interpolation operator
$\bm{P}_l$. In AMG, we take a simple smoothing scheme like Jacobi or
Gauss-Seidel relaxation. For symmetric positive definite M-matrices,
an investigation of the smooth error components
$\bm e = (e_1,\ldots,e_n)^{\top}$ reveals that these
satisfy $\bm{D}^{-1}\bm{\Lambda} \bm e \approx 0$, where $\bm D:=\diag{\bm \Lambda}$.
If we denote the coefficients of $\bm \Lambda$ by
$\lambda_{ij}$, this means
\begin{equation}
  e_i \approx - \frac{1}{\lambda_{ii}} \sum_{j \neq i} \lambda_{ij} e_j. \label{eq:smooth_error_i}
\end{equation}
Equation \eqref{eq:smooth_error_i} already delivers us a template for
the interpolation formula: The value at the fine grid point $i$ should be
approximated by the values at those points $j$ for which $- \lambda_{ij}$ is
relatively large, while those (scaled) entries $\lambda_{ij}$ also provide
the interpolation weights.  In consequence, a substantial amount of
those large negative connections should lead (directly or indirectly)
to coarse grid points $j \in \bm{C}_l$. This imposes certain conditions on
the selection of the coarse grid points $\bm{C}_l \subset \Omega_l$.
Many algorithms to the coarse grid selection and construction of the
interpolation have been proposed. We do not
describe them in detail here, but refer to \cite{RugeStueben87,
  Stueben99, HensonYang02, DeSterckYangHeys06,
  DeSterckFalgoutNoltingYang08}.
In our experiments we use Falgout coarsening to select the set of coarse grid
points $\mathfrak{C}_l$.
We compute the interpolation weights $\bm{P}_{l,\mathfrak{F}}$ using the
modified classical interpolation scheme and dropped interpolation weights
smaller than 1/20 of the largest absolute weight per row.
The restriction matrix is taken as the transpose of the interpolation.

The Galerkin ansatz for the coarse grid matrix \eqnref{galerkinproduct}
has two benefits:
First, for symmetric positive definite
$\bm{\Lambda}_l$, $\bm{\Lambda}_{l+1}$ is also symmetric positive definite for any full
(column) rank interpolation operator $\bm{P}_l$. Second, the resulting
two-grid correction operator $\bm I - \bm{P}_l \bm{\Lambda}_{l+1}^{-1} \bm{P}_l \bm{\Lambda}_l$ is an
orthogonal projector, which (extended recursively over all levels and
combined with smoothing) ensures that the multigrid cycle converges \cite{Stueben99}.

%\bram{Feel free to add to this section.  If you would like to add to other
%sections, that's great; please put it under dedicated subsections and we will
%tie it together later.}

We now present an algebraic multigrid approach for a monolithic
solution of saddle point problems of the form
\begin{equation}
    \begin{bmatrix} \bm{A} & \bm{B}^{\top} \\ \bm{B} &
    - \bm{C}\end{bmatrix}
    \begin{bmatrix}
       \bm u \\
       p
    \end{bmatrix} =
    \begin{bmatrix}
      \bm v \\
      q
    \end{bmatrix}
    .
\end{equation}
We will construct a multigrid
hierarchy of saddle point matrices indexed by $l$,
\begin{equation}
  \mathcal{A}_l = \begin{bmatrix} \bm{A}_l & \bm{B}_l^{\top} \\ \bm{B}_l &
    - \bm{C}_l\end{bmatrix},
\end{equation}
as well as block interpolation matrices $\mathcal{P}_l$.
Note that we include a lower
right block $\bm{C}_l$, which will be zero on the
finest level, $\bm{C}_1 = 0$, and non-zero for $l \geq 2$.

\subsection{Smoothers for saddle point problems}

Classical relaxation schemes like the Jacobi or Gauss-Seidel iteration
are not suitable for saddle point systems,
since the smoothing properties for
these schemes rely on the positive definiteness of the matrix.
% , which is missing here.
We present two dedicated smoothing schemes for saddle point systems
\cite{SchoeberlZulehner03}: First, a predictor-corrector algorithm, which combines
segregated sweeps over the physical components, and second, a box
relaxation scheme, where small saddle point subsystems are solved within a
global Schwarz method.

\paragraph{Uzawa relaxation.}

Our first option is a symmetric inexact Uzawa relaxation scheme. Within each iteration, a predictor $\bm{u}^{\ast}$ for the flux is
computed first, which is used to update $\bm u$. Finally, employing the updated
values of $p_l^{it+1}$, the next iterate $\bm{u}_l^{it+1}$ is computed,
\begin{subequations}
\begin{align}
  \bm{u}_l^{\ast} &= \bm{u}^{it} + \bm{\hat A}_l^{-1} \left( \bm{v} - {\bm
      A}_l \bm{u}_l^{it} - {\bm B}_l^{\top} p_l^{it}
  \right), \label{eq:uzawa_1}\\
  p_l^{it+1} &= p_l^{it} + \bm{\hat S}_l^{-1} \left( {\bm B}_l
    \bm{u}_l^{\ast} - {\bm C}_l p_l^{it} - q \right),\\
  \bm{u}_l^{it+1} &= \bm{u}_l^{it} + \bm{\hat A}_l^{-1} \left( \bm{v} - {\bm
      A}_l \bm{u}_l^{it} - {\bm B}_l^{\top} p_l^{it+1} \label{eq:uzawa_3}
  \right).
\end{align}
\label{eq:uzawa}
\end{subequations}
The matrices $\bm{\hat A}_l$ and $\bm{\hat S}_l$ are chosen such that
$\bm{\hat A}_l - {\bm A}_l$ and $\bm{\hat S}_l - {\bm B}_l
\bm{\hat{A}}_l^{-1} {\bm B}_l^{\top} - {\bm C}_l$ are
symmetric positive definite. Furthermore, they should be easily
invertible. For example, we can choose to use the scaled diagonals of $\bm{A}_l$ and $\bm{B}_l
\bm{\hat A}_l^{-1} \bm{B}_l^{\top} + \bm{C}_l$, respectively. The magnitude of the scaling
can be obtained with help of the power iteration on the latter matrices.
%Largest eigenvalue of A and \hat{A}-A comparable !
For convergence and further properties of this smoother we refer to
\cite{SchoeberlZulehner03}.

\paragraph{Vanka smoothing.}

The second alternative requires us to first decompose the
computational domain into small overlapping patches. To this end, we
employ the non-zero structure of $\bm{B}_l$ and construct a patch $\Omega_{l,j}$
for each row of $\bm{B}_l$: For the $j$-th row, $\Omega_{l,j}$ consists of
the index $j$ as well as all indices $i$ such that there exists an
entry $b_{ji} \neq 0$,
%\todo{Show an image to demonstrate what this means for this discretization}
\begin{equation}
  \Omega_{l,j} := \{ i: \  b_{ji} \neq 0 \} \times \{ j \}.
\end{equation}
The transfer between the global domain $\Omega_l$ and the subdomains
$\Omega_{l,j}$ is accomplished by (optionally scaled) injection
operators $\bm{V}_{l,j}$ (flux) and $\bm{W}_{l,j}$ (pressure),
\begin{subequations}
\begin{align}
  \bm{V}_{l,j} &= \diag{\mathrm{\bf v}_{l,i}}_{i=1,\ldots,n_l} \bm
                 J_j \label{eq:zulehner_def_V}, \\
  \bm{W}_{l,j} &= \bm{\tilde{J}}_j  \label{eq:zulehner_def_W},
\end{align}
\end{subequations}
where $\bm J_j$ and  $\bm{\tilde{J}}_j$  are binary matrices that map the subdomain into the
global domain for velocity and pressure, respectively; \ie each of their
columns contains exactly one unit entry.
%zero and are zero otherwise.

First, the residuals are restricted to the subdomain,
\begin{subequations}
\begin{align}
  \bm{v}_{l,j} &= \bm{V}^{\top}_{l,j} \left( \bm{v} - \bm{A}_l
                 \bm{u}_l - \bm{B}_l^{\top} p_l \right), \label{eq:vanka_1}\\
  q_{l,j} &= \bm{W}^{\top}_{l,j} \left( q - \bm{B}_l \bm{u}_l + \bm{C}_l^{\top} p_l \right).
\end{align}
\end{subequations}
Then, on each subdomain, a small saddle point problem of the form
\begin{equation}
    \begin{bmatrix} \bm{\hat A}_{l,j} & \bm{B}_{l,j}^{\top} \\ \bm{B}_{l,j} & \bm{B}_{l,j} \bm{\hat A}_{l,j}^{-1} \bm{B}_{l,j}^{\top}
    - \bm{\hat S}_{l,j} \end{bmatrix} \begin{bmatrix} \bm{u}_{l,j} \\ p_{l,j} \end{bmatrix}
  = \begin{bmatrix} \bm{v}_{l,j} \\ q_{l,j} \end{bmatrix}
\end{equation}
is solved using a sparse direct solver.
Finally, the updates are prolongated to the global domain,
\begin{subequations}
\begin{align}
  \bm{u}_l &= \bm{u}_l + \bm{V}_{l,j} \bm{u}_{l,j},\\
  p_l &= p_l + \bm{W}_{l,j} p_{l,j}. \label{eq:vanka_5}
\end{align}
\end{subequations}
The small matrices $\bm{\hat A}_{l,j}$, $\bm{B}_{l,j}$ and
$\bm{\hat S}_{l,j}$, where the latter is just a scalar, are defined by
\begin{subequations}
\begin{align}
  \bm{\hat A}_{l,j} &=   \bm{V}_{l,j}^{\top} \bm{\hat A}_l\bm{\hat V}_{l,j}, \label{eq:zulehner_def_Ahatj}\\
  \bm{B}_{l,j} &= \bm{W}_{l,j}^{\top} \bm{B}_l \bm{\hat V}_{l,j}, \label{eq:zulehner_def_Bj}\\
  \bm{\hat S}_{l,j} &= \beta^{-1} ( \bm{C}_{l,j} + \bm{B}_{l,j} \bm{\hat
A}_{l,j}^{-1} \bm{B}_{l,j}^{\top} ), \label{eq:zulehner_def_Shatj}
\end{align}
where $\bm{C}_{l,j}$ denotes the $j$-th diagonal entry of $\bm{C}_l$
and $\bm{\hat V}_{l,j} = \diag{\mathrm{\bf v}_{l,i}^{-1}}_{i=1,\ldots,n_l} \bm J_j$,
the inverse scaling of \eqref{eq:zulehner_def_V}.
\end{subequations}
Again, $\bm{\hat A}_l$ is a scaled version of the diagonal of ${\bm
 A}_l$ such that $\bm{\hat A}_l - {\bm A}_l$ is positive
definite.
With the aid of a power iteration, $\beta > 0$ is computed such that
\begin{equation}
  \diag{\bm{\hat S}_{l,j}}_i  - \left( \bm{C}_l + \bm{B}_l
  \bm{\hat A}_l^{-1} \bm{B}_l^{\top} \right)
\end{equation}
is symmetric positive definite.

The iteration can be performed either additively (\ie the residuals
are computed once, then all subdomain solves are performed
independently) or multiplicatively (after each subdomain solve the
residuals are updated). In the latter case, it is beneficial to
perform a symmetric sweep, \ie after one complete sweep the
subdomain solves are performed again in reverse order.  If we choose
\begin{align}
  \mathrm{\bf v}_{l,i} = \frac{1}{\sqrt {\vert \{ j: b_{ji} \neq 0 \} \vert}},\label{eq:sz_scale}
\end{align}
the additive smoother coincides with the Uzawa method described in the
previous section \cite{SchoeberlZulehner03}. On the other hand, in the
multiplicative case the simple choice
\begin{align}
\mathrm{\bf v}_{l,i} = 1
\label{eq:sz_one}
\end{align}
may result in faster convergence.

\subsection{AMG setup for saddle point systems}

The starting point for our saddle point AMG method (SPAMG) is the block
diagonal matrix
\begin{equation}
	\mathcal{B}=
  \begin{bmatrix}
    \bm{A}_l & 0 \\ 0 & \bm{B}_l \bm{\hat A}_l^{-1} \bm{B}_l^{\top} + \bm{C}_l
\end{bmatrix}.
\end{equation}
This operator is symmetric positive definite.
%\bram{When does
%  the ``semi'' apply? BM: one of the block diagonals may have a kernel.
% Should we really discuss that here?}.
Let us assume
that we can apply the classical AMG setup algorithm as described in
Section \ref{section:amg_elliptic} to each of the blocks $\bm{A}_l$ and
$\bm{B}_l \bm{\hat A}_l^{-1} \bm{B}_l^{\top} + \bm{C}_l$.
We construct coarse grids and interpolation operators $\bm{P}_{l,u}$ and
$\bm{P}_{l,p}$ for each of these blocks and obtain a block interpolation
operator
\begin{equation}
\mathcal{P}_l :=
  \begin{bmatrix}
    \bm{P}_{l,u} & 0 \\ 0 & \bm{P}_{l,p}
\end{bmatrix}.
\end{equation}
Now, one idea would be to construct the coarse grid operator in
the usual Galerkin way, $\mathcal{A}_{l+1} := \mathcal{P}^{\top}_l
\mathcal{A}_l \mathcal{P}_l$.  Unlike in the symmetric positive definite
case, however, we cannot be sure that $\mathcal{A}_{l+1}$ is invertible,
thus we might attain an unusable multigrid hierarchy.
To prevent this, we
modify the interpolation operator $\mathcal{P}_l$ by the
multiplication of a stabilization term. To this end, we re-write the
velocity prolongation block such that the rows corresponding to the
fine grid points come first,
\begin{equation}
 \bm{P}_{l,u} =  \begin{bmatrix} \bm{P}_{l,u,CF} \\
    \bm{I}_{l,u,C} \end{bmatrix},
\end{equation}
where $\bm{I}_{l,u,C}$ is the identity injecting the values from
level $l+1$ to level $l$ and $\bm{P}_{l,u,CF}$ contains the interpolation
weights from coarse to fine computed from the matrix $\bm A_l$. Analogously, we write
\begin{equation}
  \bm{B}^{\top}_l = \begin{bmatrix} \bm{B}^{\top}_{l,F} \\ \bm{B}^{\top}_{l,C} \end{bmatrix}, \quad
  \bm{\hat A}^{\top}_l = \begin{bmatrix} \bm{\hat A}^{\top}_{l,F} & 0 \\ 0 & \bm{\hat A}^{\top}_{l,C} \end{bmatrix}.
\end{equation}
The stabilized interpolation operator is computed by
\begin{equation}
  \tilde{\mathcal{P}}_l := \begin{bmatrix} \bm{I}_{l,F} & 0 & - \bm{\hat A}_{l,F}^{-1}
    \bm{B}_{l,F}^{\top} \\ 0 & \bm{I}_{l,C} & 0 \\
    0 & 0 & \bm{I}
  \end{bmatrix}
  \begin{bmatrix}
    \bm{P}_{l,u,CF} & 0 \\
    \bm{I}_{l,u,C} & 0 \\
    0 & \bm{P}_{l,p}
  \end{bmatrix}.
\end{equation}
% Interpolation conditions on the pressure are weaker, so we don't need to split it like the velocity.
Now, we use the modified Galerkin product
\begin{equation}
  \eqnlab{modgalerkin}
  \mathcal{A}_{l+1} : = \tilde{\mathcal{P}}^{\top}_l
  \mathcal{A}_l \tilde{\mathcal{P}}_l
\end{equation}
to compute the coarse grid
operator. For this matrix, we can show an inf-sup-condition if (a) the
fine grid matrix $\mathcal{A}_l$ fulfils such a condition and (b) the
interpolation operator $\bm{P}_{l,u}$ for the velocity satisfies certain
approximation properties, which most usual AMG interpolation schemes
do \cite[Lemma 4.6]{Metsch13}.
The invertibility of the coarse grid matrix \eqnref{modgalerkin} is ensured by
its lower right block
% $-\bm{C}_{l+1}$,
\begin{equation}
  {\bm C}_{l+1} = \bm{P}^{\top}_{l,p} \left[ \bm{C}_l + 2 \bm{B}_{l,F} \bm{\hat A}_{l,F}^{-1} \bm{B}^{\top}_{l,F}+
    \bm{B}_{l,F} \bm{\hat A}_{l,F}^{-1} \bm{A}_{l,F} \bm{\hat A}_{l,F}^{-1} \bm{B}_{l,F}^{\top} \right] \bm{P}_{l,p}.
\end{equation}
This can be understood as a partial Schur complement that ensures the stability of the coarse grid matrix.

In reference to recent work, the results published in
\cite{Webster:2015} and \cite{Webster:2018} indicate that it is
neither always necessary to stabilize the interpolation operator on every
level, nor that interpolation and restriction need to be stabilized
simultaneously.
A further investigation of these techniques, however, is beyond the scope of
this paper:
It may be quite involved to determine algorithmically and separately for each
level whether the stabilization is required.

\section{Numerical Results}
\seclab{numerical}

% \todo{CB: Work into the text below:
%       The focus of our paper is neither error estimation nor AMR.
%       We use AMR purely to construct a more difficult system to
%       solve by the various preconditioners we compare.}

In this section, we evaluate the effectiveness of SPAMG compared to the diagonal
and Schur preconditioners for uniform and adaptive meshes. We choose the examples
based on manufactured solutions, that is, we prescribe a target pressure field and
compute the velocity field, boundary conditions, and right hand sides in order
to satisfy \eqnref{poisson:system}.
We also choose synthetic adaptive refinement tailored to the reference solutions.
Examples 1 to 3 are typical benchmarks with smooth coefficients.
In example 4 we challenge the
numerical solver and the preconditioner by using a conductivity tensor with strong
coefficient variation. We employ three different flavors of smoothers
inside SPAMG: The inexact Uzawa scheme
\eqref{eq:uzawa} and two variants of the
multiplicative Vanka smoother, either using the scaling
\eqref{eq:sz_scale} (``Vanka Scale") or the unscaled injection
\eqref{eq:sz_one} (``Vanka One").
%\jose{You always employed the Vanka smoother multiplicatively, didn't
%  you? JF: Yes}
%\todo{CB/JF: Is Vanka one / Vanka scale sufficiently explained?}

We delegate the parallel mesh management
to the octree-based AMR software library \pforest
\cite{BursteddeWilcoxGhattas11, IsaacBursteddeWilcoxEtAl15}. This library
provides a collection of algorithms that implement scalable parallel AMR
operations. In particular, we employ \pforest  to create and modify
a hexahedral triangulation of the unit square/cube and to introduce a
numbering of the degrees of freedom suitable for a \rtzero discretization. Regarding
solvers and preconditioners, we use the GMRES \cite{SaadSchultz86a}
implementation provided by the software library \hypre \cite{HypreTeam12}.
In order to approximate the product $\bm N \bm r$ for a given
residual $\bm r$, we employ one SPAMG V-cycle.
We use the parallel multigrid implementation BoomerAMG from \hypre as the frame
to which our saddle-point AMG is added.
% Finally, the SPAMG preconditioner
% was implemented inside a private branch of \hypre.

% We  visualize the pressure and velocity fields arising from
% the \rtzero discretization using the VTK file format. The pressure scalar field,
% being approximated by discontinuous piecewise constants is trivially
% written as cell data. For the velocity vector field, we require to save
% it as point data. The VTK routines require information for each element corner,
% while there are $2^d$ corners per element, we have $2$ velocity degrees
% of freedom per coordinate direction for $d=2,3$.
% We have addressed this mismatch
% by replicating the computed velocity
% value to the corners of the face containing it. In \figref{rt0_vis}
% we show the result for a two dimensional vector field.
% \begin{figure}
% 	\centering
% 	\subfloat[$x$ components]{\includegraphics[scale=1.5]{rt0_x_vis.pdf}}
% 	\hspace{3em}
% 	\subfloat[$y$ components]{\includegraphics[scale=1.5]{rt0_y_vis.pdf}}
% 	\caption{\rtzero dof identification for VTK visualization of a two
% 		dimensional vector field. Bold curved arrows indicate value copy.}
% 	\figlab{rt0_vis}
% \end{figure}

\subsection{Homogeneous Dirichlet boundary conditions}
\seclab{example1}
%\jose{Change these subsection titles into something short, meaningful, and
%      avoid repetitions.}
%\jose{Use proper sentences, but keep it short.}
We solve \eqnref{poisson:system} on the unit cube using an identity
conductivity tensor and homogeneous Dirichlet boundary conditions.
We compute the right hand side based on the manufactured solution
\begin{equation}
	p(x,y) = (x^2 - x^3)(y^2 - y^3)
\end{equation}
in 2D and
\begin{equation}
	p(x,y,z) = (x^2 - x^3)(y^2 - y^3)(z-z^2)
\end{equation}
in 3D, respectively.
\begin{figure}
	\centering
	\subfloat[]{\includegraphics[scale=0.85]{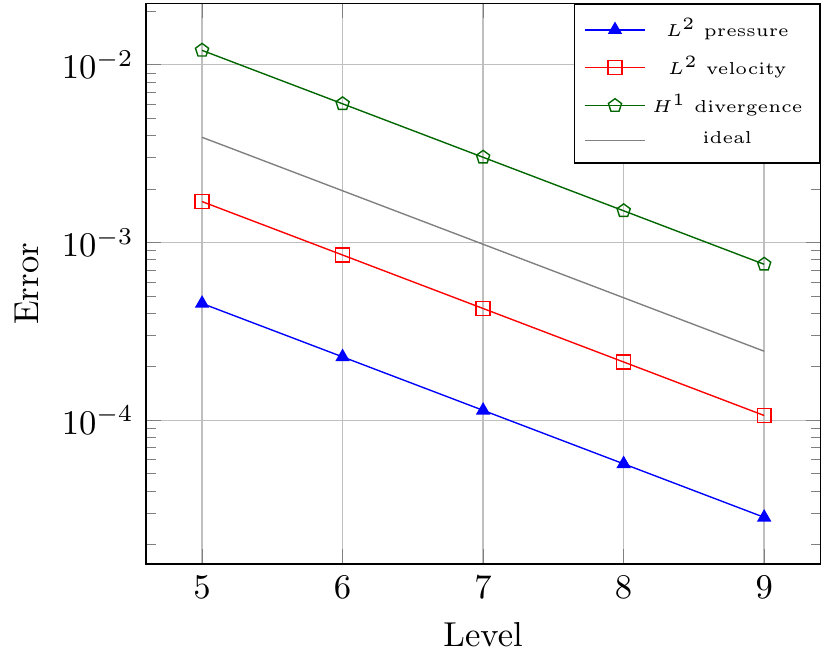}}
	\hfill
	\subfloat[]{\includegraphics[scale=0.85]{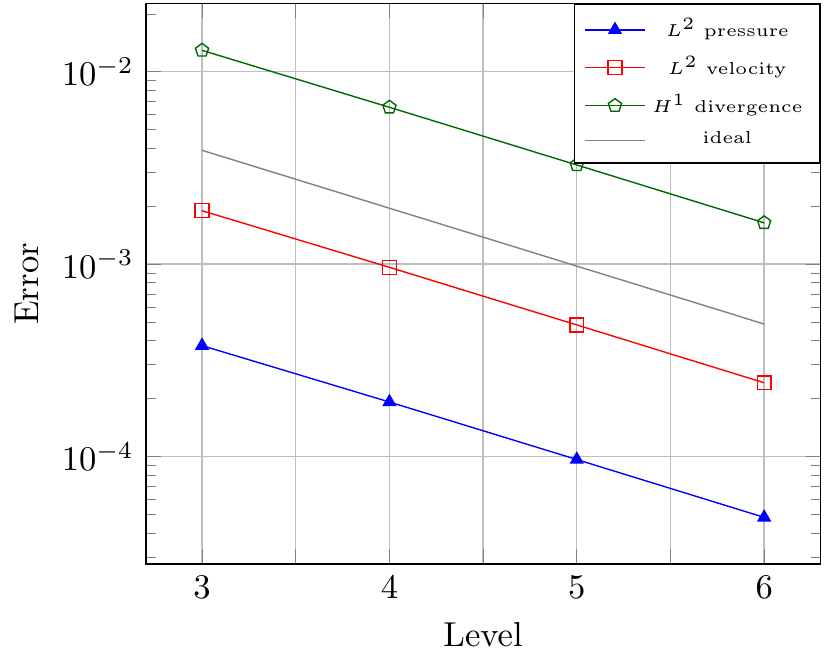}}
	\caption{Error plot for the numerical solution of a mixed Poisson system
		corresponding to the example specified in \secref{example1}
		(homogeneous Dirichlet boundary conditions) in two (a)
		and three dimensions (b) for uniform meshes and identity conductivity tensor.
		The level $\ell$ is related to the mesh size $h$ via $h = 2^{-\ell}$.
		We confirm the expected
		convergence rates predicted by \eqnref{conv_rate}.}
	\figlab{err_e1_u}
\end{figure}
\begin{figure}
	\centering
	\subfloat[]{\includegraphics[scale=0.85]{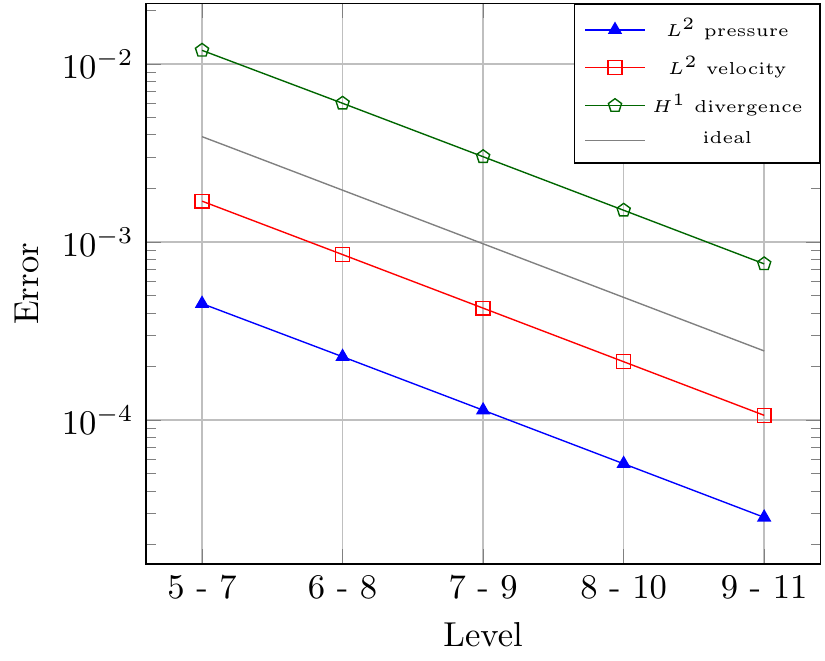}}
	\hfill
	\subfloat[]{\includegraphics[scale=0.85]{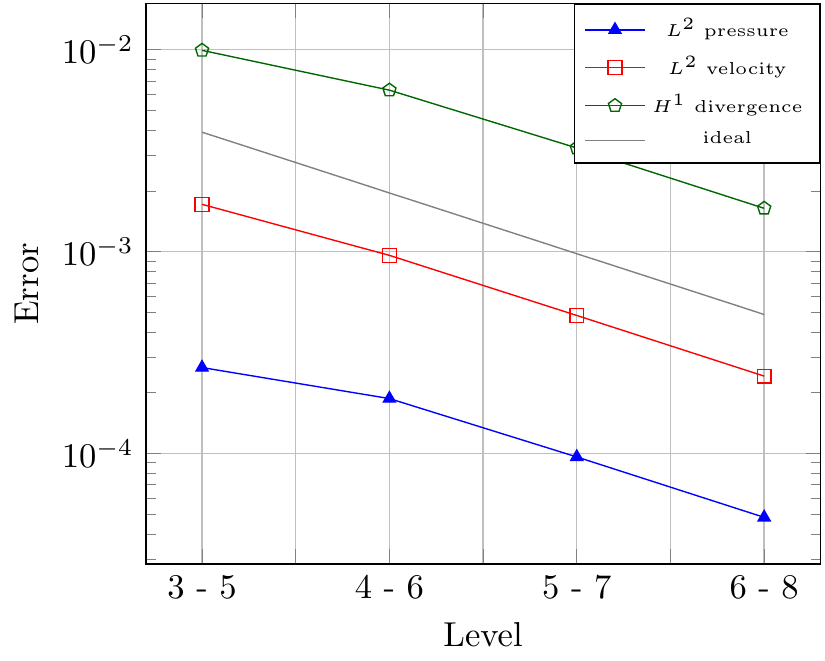}}
	\caption{Error plot
                % for the numerical solution of a mixed Poisson system
		% corresponding to the example specified in \secref{example1}
                in two (a) and three dimensions (b) for adaptive meshes and
                identity conductivity tensor.
		We choose to refine an element
		of side length $h$ by two additional levels
		whenever its centroid lies within the circle/sphere of
		radius $2^{d}h^{2}$ around the point $(\frac{1}{2},\frac{1}{2})$ for $d=2$
		and $(\frac{1}{2},\frac{1}{2}, \frac{1}{2})$ for $d=3$.
		Because of the smoothness of the solution and the
		equations coefficients we do not expect that local refinement
		translates into an improvement of the approximation with respect to a
		uniform case mesh.}
	\figlab{err_e1_a}
\end{figure}
\begin{table}
	\renewcommand{\tabcolsep}{1ex}
	\renewcommand{\arraystretch}{1.2}
	\centering
	\subfloat[]{\scalebox{0.65}{%
	\begin{tabular}{ ccccccc }
		\hline
		Level	& \multicolumn{6}{c}{\# Iterations} \\ \cline{2-7}
			& NoPC 	& Diag 	& Schur & \multicolumn{3}{c}{SPAMG}  \\
			&	&	&	& Uzawa	& Vanka one	& Vanka scale \\ \cline{5-7}
		$4$	& 459	& 299  	& 22  	& 9	& 8	& 8	\\
		$5$	&$>$1000& 773  	& 22	& 10	& 8	& 8 	\\
		$6$	& -	&$>$1000& 24	& 12	& 9	& 10	\\
		$7$	& -	& -	& 24	& 14	& 10	& 10	\\
		$8$	& -	& -	& 24	& 14	& 10	& 11	\\
		$9$	& -	& -	& 24	& 14	& 10	& 11	\\
		\hline
	\end{tabular}
	}}
	\hfill
	\subfloat[]{\scalebox{0.65}{%
	\begin{tabular}{ ccccccc }
		\hline
		Level & \multicolumn{6}{c}{\# Iterations} \\ \cline{2-7}
			& NoPC 	& Diag 	& Schur & \multicolumn{3}{c}{SPAMG}  \\
			&	&	&	& Uzawa	& Vanka one	& Vanka scale  \\ \cline{5-7}
		$4-7$   &$>$1000&$>$1000& 87	& 10	 & 7	&	7	\\
		$5-8$   & -	& -  	& 87	& 10	 & 7	&	8 	\\
		$6-9$   & -	& - 	& 99	& 13	 & 10	&	10	\\
		$7-10$  & -	& -	& 98	& 13	 & 10	&	10	\\
		$8-11$  & -	& -	& 112	& 15	 & 11	&	11	\\
		$9-12$	& -	& -	& 149	& 15	 & 11	&	11	\\
		\hline
	\end{tabular}
	}}
	\caption{Number of iterations required by the GMRES solver for a two
		dimensional mixed Poisson system (\secref{example1})
		on a uniform (a) and adaptive mesh (b).
		The relative tolerance of the linear solve is fixed to  $10^{-6}$.
		We display iteration counts employing no preconditioning in column two,
                a diagonal lumped mass matrix in column three, a Schur
                complement preconditioner
		in column four and three different smoothers in the
		SPAMG preconditioner in columns five to seven from each table.}
	\tablab{ex1_iter_2d}
\end{table}
\begin{table}
	\renewcommand{\tabcolsep}{1ex}
	\renewcommand{\arraystretch}{1.2}
	\centering
	\subfloat[]{\scalebox{0.85}{%
	\begin{tabular}{ ccccc }
		\hline
		Level & \multicolumn{4}{c}{\# Iterations} \\ \cline{2-5}
				& Schur & \multicolumn{3}{c}{SPAMG}  \\
				&	& Uzawa	& Vanka one   & Vanka scale \\
		$3$   &	20	& 9	& 7	&	7	\\
		$4$   &	22	& 11	& 8	&	8	\\
		$5$   &	22 	& 13	& 9	&	9 	\\
		$6$   & 24  	& 14	& 9	&	10	\\
		$7$   & 24  	& 15	& 11	&	12 	\\
		\hline
	\end{tabular}
	}}
	\hfill
	\subfloat[]{\scalebox{0.85}{%
	\begin{tabular}{ ccccc }
		\hline
		Level & \multicolumn{4}{c}{\# Iterations} \\ \cline{2-5}
			& Schur & \multicolumn{3}{c}{SPAMG}  \\
			& 	& Uzawa  & Vanka one   & Vanka scale \\
		$3-6$   & 238	& 11	& 8	& 8	\\
		$4-7$   & 322	& 11	& 8	& 8 	\\
		$5-8$   & 338	& 13	& 9	& 9	\\
		$6-9$   & 383	& 14	& 10	& 10	\\
		$7-10$  &   390	& 17	& 12	& 12 	\\
		\hline
	\end{tabular}
	}}
	\caption{Number of iterations required by the GMRES solver for a three
		dimensional mixed Poisson system (\secref{example1}) discretized on a
		uniform (a) and adaptive mesh (b).% We use the same setup as in
		%\tabref{ex1_iter_2d}.
                }
	\tablab{ex1_iter_3d}
\end{table}

%\jose{Very short summary of main features we see in the results.
%      Do this for all examples.}
With this example we verify the correctness of our implementation of the $\rtzero$ discretization.
The discretization error converges the predicted rates as confirmed in \figref{err_e1_u} and \figref{err_e1_a}.
The iteration counts displayed in \tabref{ex1_iter_2d} and \tabref{ex1_iter_3d} confirm the
(well known) robustness of the Schur preconditioner for uniform meshes. For adaptive meshes,
the Schur complement preconditioner incurs high iteration counts, particularly for the 3d case.
The three variants of the SPAMG preconditioner retain
mesh independent iteration counts for both uniform and adaptive meshes.

\subsection{Inhomogeneous Dirichlet/Neumann boundary conditions}
\seclab{example2}
We solve \eqnref{poisson:system} with an identity tensor. We impose homogeneous
mixed homogeneous Dirichlet / Neumann boundary conditions and compute the right
hand side based on the exact solution
\begin{equation}
p(x,y) = xy(1-y)(1-x)^2
\end{equation}
in 2D and
\begin{equation}
p(x,y,z) = xy(1-y)(1-x)^2(1 - z)
\end{equation}
in 3D, respectively. The Neumann boundary is set at $y=0$ and $y=1$  in both cases.
\begin{figure}
	\centering
	\subfloat[]{\includegraphics[scale=0.85]{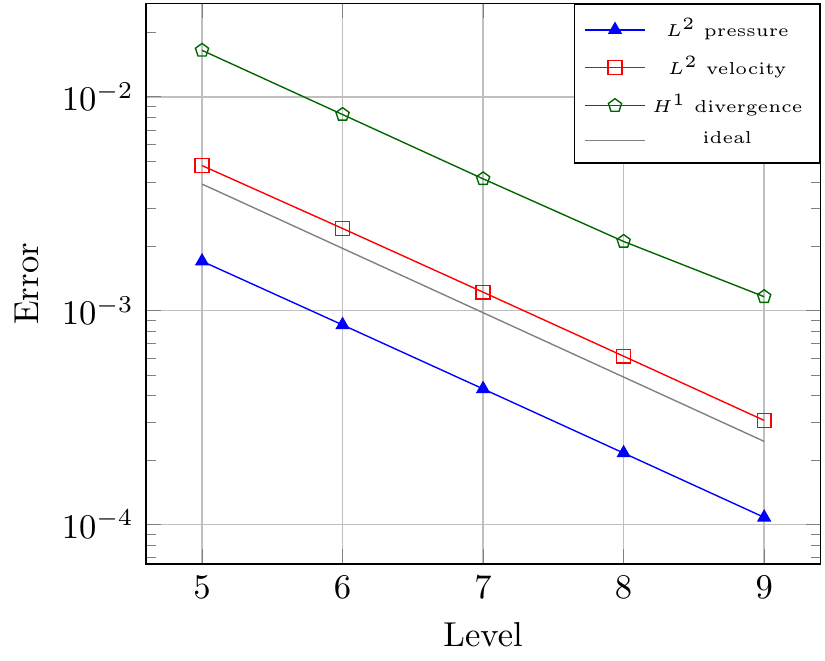}}
	\hfill
	\subfloat[]{\includegraphics[scale=0.85]{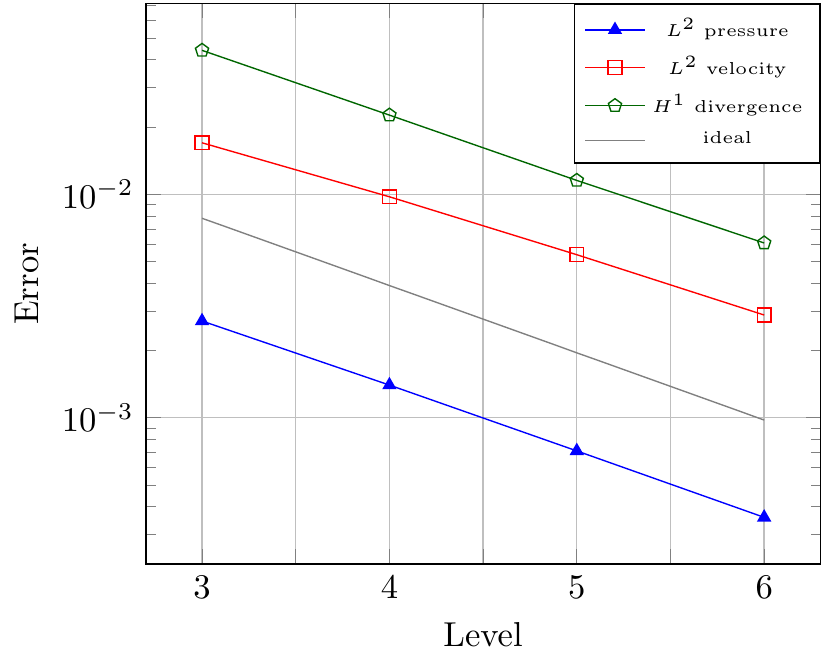}}
	\caption{Error plot
		(\secref{example2}, inhomogeneous Dirichlet/Neumann boundary conditions)
		in two (a) and three (b) dimensions for uniform meshes.}
	\figlab{err_e2_u}
\end{figure}
\begin{figure}
	\centering
	\subfloat[]{\includegraphics[scale=0.85]{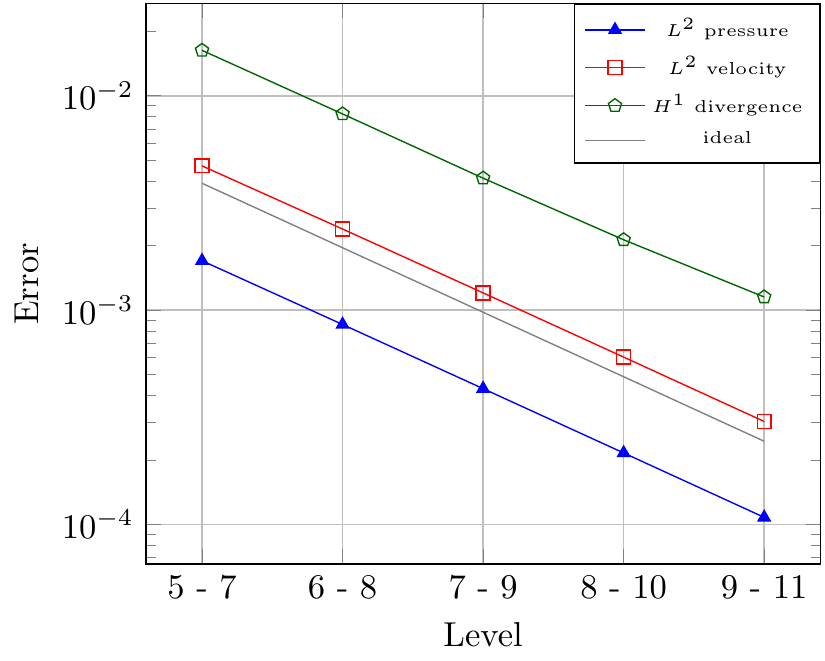}}
	\hfill
	\subfloat[]{\includegraphics[scale=0.85]{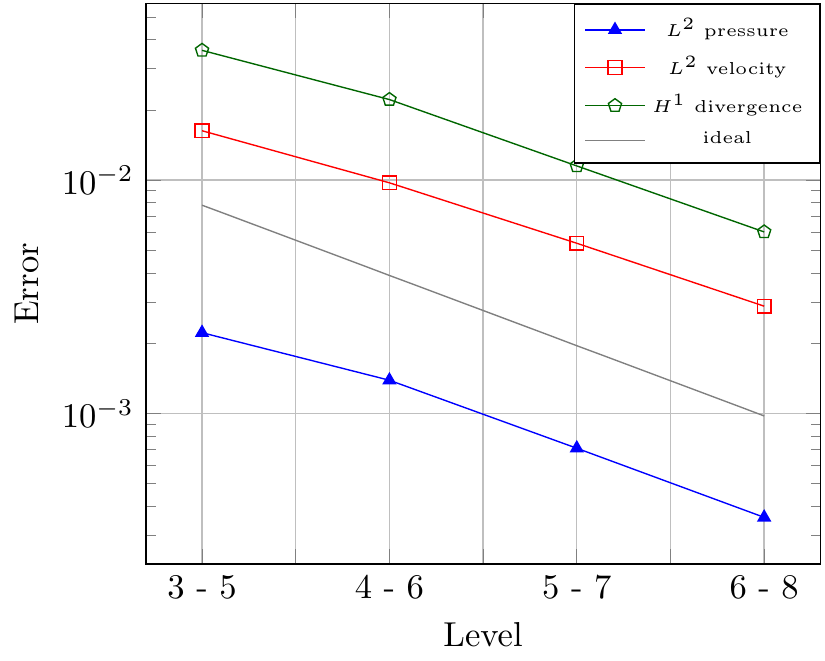}}
	\caption{Error plot
        % for the numerical solution of a mixed Poisson system
	% specified in
        (\secref{example2})
	in two (a) and three (b) dimensions for adaptive meshes.
	The refinement criterion is chosen as in \figref{err_e1_a}.}
	\figlab{err_e2_a}
\end{figure}
%In this example we check the correctness of our implementation
%of the approach \eqnref{neumann_enforce} to enforce the Neumann boundary
%condition \eqnref{poisson:system:neumann}.
As in the previous section, the theoretical converge rates agree with the bounds
\eqnref{conv_rate}. Due to the smoothness solution and the
equation constant coefficient, no additional benefit is expected from local adaptation of the mesh; see \figref{err_e2_u} and \figref{err_e2_a}.
The iteration counts (not shown) are similar to the previous example.

\subsection{A non-trivial conductivity tensor}
\seclab{example3}
In this example we approximate the solution of \eqnref{poisson:system} for the case of a non-diagonal
conductivity tensor $\K$. We impose non-zero Dirichlet boundary conditions.\
The manufactured solution is
\begin{equation}
	p(x,y) = e^x\sin(y)
\end{equation}
in 2D and
\begin{equation}
	p(x,y,z) = e^x\sin(y)(1+z^2)
\end{equation}
in 3D, respectively. The conductivity tensor is given by
\begin{equation}
	\K(x,y) =
	\begin{pmatrix}
	e^{x/2 + y/4} & \sin(2\pi x) \\
	\sin(2\pi x)  & e^{x/4 + y/2}
	\end{pmatrix}
\end{equation}
in 2D and
\begin{equation}
	\K(x,y,z) =
	\begin{pmatrix}
	e^{x/2 + y/4} & \sin(2\pi x) & 0 \\
	\sin(2\pi x)  & e^{x/4 + y/2} & 0 \\
	0 & 0 & e^{z}
	\end{pmatrix}
\end{equation}
in 3D.
% \begin{figure}
% 	\centering
% 	\subfloat[]{\includegraphics[scale=0.85]{e3_2d_unif.pdf}}
% 	\hfill
% 	\subfloat[]{\includegraphics[scale=0.85]{e3_3d_unif.pdf}}
% 	\caption{Error plot for the numerical solution of a mixed Poisson system
% 	specified in \secref{example3} in two (a) and three (b) dimensions for uniform meshes.}
% 	\label{err_e3_u}
% \end{figure}
% \begin{figure}
% 	\centering
% 	\subfloat[]{\includegraphics[scale=0.85]{e3_2d_adapt.pdf}}
% 	\hfill
% 	\subfloat[]{\includegraphics[scale=0.85]{e3_3d_adapt.pdf}}
% 	\caption{Error plot for the numerical solution of a mixed Poisson system
% 	specified in \secref{example3}  in two (a) and three (b) dimensions for adaptive meshes.
% 	The refinement criteria is chosen as in \figref{err_e1_a}.}
% 	\label{err_e3_a}
% \end{figure}
\begin{table}
	\renewcommand{\tabcolsep}{1ex}
	\renewcommand{\arraystretch}{1.2}
	\centering
	\subfloat[]{\scalebox{0.65}{%
	\begin{tabular}{ ccccccc }
		\hline
		Level & \multicolumn{6}{c}{\# Iterations} \\ \cline{2-7}
		& NoPC 	& Diag 	& Schur & \multicolumn{3}{c}{SPAMG}  \\
		&	&	& 	& Uzawa	& Vanka one	& Vanka scale \\ \cline{5-7}
		$4$	&$>$1000&$>$1000& 85	& 18	& 12	& 11	\\
		$5$	& -	& -  	& 102	& 20	& 13	& 13	\\
		$6$	& - 	& -	& 113	& 23	& 14    & 14	\\
		$7$	& - 	& -	& 130	& 24	& 16    & 15	\\
		$8$	& - 	& -	& 142	& 25	& 16    & 17	\\
		$9$	& - 	& -	& 169	& 25	& 17    & 18	\\
		\hline
	\end{tabular}
	}}
	\hfill
	\subfloat[]{\scalebox{0.65}{%
	\begin{tabular}{ ccccccc }
		\hline
		Level & \multicolumn{6}{c}{\# Iterations} \\ \cline{2-7}
		& NoPC 	& Diag 	& Schur & \multicolumn{3}{c}{SPAMG}  \\
		&	&	& 	& Uzawa	& Vanka one	& Vanka scale  \\ \cline{5-7}
		$4-7$   &$>$1000&$>$1000& 250 	& 18	& 11	& 10	\\
		$5-8$   & -	& -  	& 255	& 19 	& 13	& 12	\\
		$6-9$   & -	& - 	& 415	& 22	& 14    & 13	\\
		$7-10$  & -	& -		& 338	& 24    & 15	& 15	\\
		$8-11$  & -	& -		& 394 	& 24	& 16    & 16	\\
		$9-12$	& -	& -		& 539 	& 24	& 16    & 17	\\
		\hline
	\end{tabular}
	}}
	\caption{Number of iterations required by the GMRES solver for a two
                dimensional mixed Poisson system with non-trivial coefficient
                tensor (\secref{example3}) discretized on
		a uniform (a) and adaptive mesh (b). We use the same setup as in
		\tabref{ex1_iter_2d}.}
	\tablab{ex3_iter_2d}
\end{table}
\begin{table}
	\renewcommand{\tabcolsep}{1ex}
	\renewcommand{\arraystretch}{1.2}
	\centering
	\subfloat[]{\scalebox{0.85}{%
	\begin{tabular}{ ccccc }
		\hline
		Level & \multicolumn{4}{c}{\# Iterations} \\ \cline{2-5}
			& Schur & \multicolumn{3}{c}{SPAMG}  \\
			&	& Uzawa	& Vanka one  	& Vanka scale \\
		$3$   &	65	& 19	& 13	& 13	\\
		$4$   &	84	& 21	& 14	& 14	\\
		$5$   &	103 & 23	& 15	& 16	\\
		$6$   & 113	& 24	& 17	& 17	\\
		$7$   & 128 & 26	& 18	& 19 	\\
		\hline
	\end{tabular}
	}}
	\hfill
	\subfloat[]{\scalebox{0.85}{%
		\begin{tabular}{ ccccc }
			\hline
			Level & \multicolumn{4}{c}{\# Iterations} \\ \cline{2-5}
				& Schur & \multicolumn{3}{c}{SPAMG}  \\
				& 	& Uzawa	& Vanka one	& Vanka scale \\
			$3-6$   & 330   & 22	& 14	& 14	\\
			$4-7$   & 456	& 23	& 15	& 14 	\\
			$5-8$   & 486	& 24	& 16	& 15	\\
			$6-9$   & 471	& 25	& 17	& 18	\\
			$7-10$  & 493   & 27	& 18	& 20	\\
			\hline
		\end{tabular}
	}}
	\caption{Number of iterations required for a three
		dimensional mixed Poisson system with non-trivial conductivity
                tensor (\secref{example3}), discretized on a
		uniform (a) and adaptive mesh (b).
                % We use the same setup as in
		%\tabref{ex1_iter_2d}.
                }
	\tablab{ex3_iter_3d}
\end{table}
%The example presented in this section aims to test the effectiveness of the SPAMG
%preconditioner for a full conductivity tensor example.
The three variants of
the SPAMG preconditioner offer mesh independent iteration counts for uniform and
adaptive meshes.
The Schur complement again produces growing iteration counts, in particular for
the three dimensional case and even for uniform meshes.
See \tabref{ex3_iter_2d} and \tabref{ex3_iter_3d}.

\subsection{High conductivity contrast}
\begin{figure}
	\centering
	\includegraphics[scale=0.9]{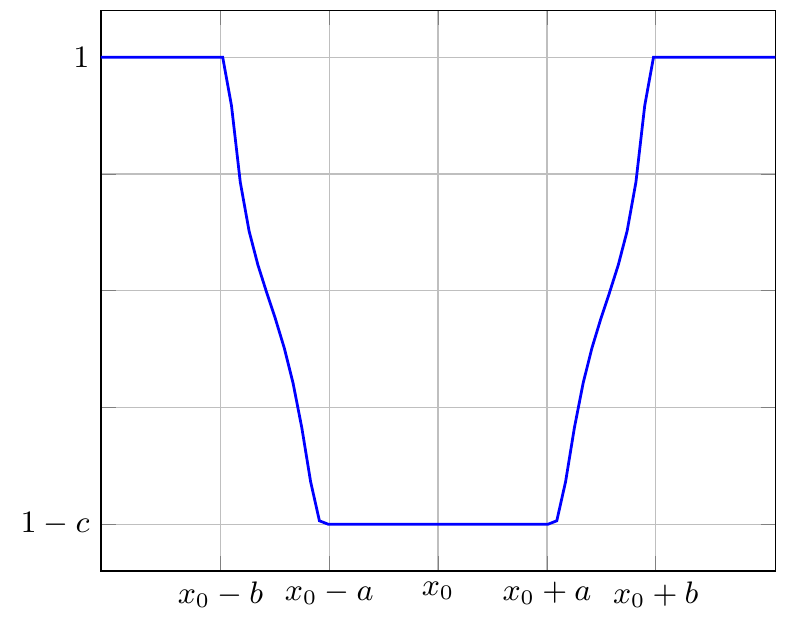}
	\caption{Sample plot of the one dimensional version of $m(\bm{x}; \bm{x_0}, a, b, c)$.}
  \label{fig:mcfunction}
\end{figure}

\seclab{example4}
To further examine the robustness of the SPAMG preconditioner,
we solve \eqnref{poisson:system} with a  conductivity tensor that exhibits strong
coefficient variation. We enforce inhomogeneous mixed Dirichlet boundary conditions
and compute the right hand side terms based on the  manufactured solution
\begin{equation}
p(x,y) = \sin(x) e^{y}
\end{equation}
in 2D and
\begin{equation}
p(x,y,z) =\sin(x) e^{y}(1+z^2)
\end{equation}
in 3D. The conductivity Tensor $\K$ is given by the identity matrix scaled
pointwise with a continuously differentiable function $m(\bm{x}; \bm{x}_0, a,
b, c)$ constructed to fulfil the following properties:
\begin{itemize}
	\item $m(\bm{x}; \bm{x}_0, a, b, c) = 1\:$ if
		$\norm{\bm{x}-\bm{x}_0}{}\geq b$,
	\item $m(\bm{x}; \bm{x}_0, a, b, c) = 1-c\:$ if
		$\norm{\bm{x}-\bm{x}_0}{} \leq a$ and
	\item $m(\bm{x}; \bm{x}_0, a, b, c) \in (1-c ,1)\:$  if
		$a < \norm{\bm{x}-\bm{x}_0}{} < b$.
\end{itemize}
Such a function can be constructed by defining
\begin{equation}
	h(t) :=
	\begin{cases}
	e^{-1/t} &\text{if } t> 0, \\
	0 &\text{else}
	\end{cases}
\end{equation}
and
\begin{equation}
	m(\bm{x}; \bm{x_0}, a, b, c) := 1 - c
	\frac{h(b-\norm{\bm{x}-\bm{x}_0}{})}
		{h(b-\norm{\bm{x}-\bm{x}_0}{}) + h(\norm{\bm{x}-\bm{x}_0}{} - a)}
	;
\end{equation}
see \figref{mcfunction}.
Hence, the first three arguments define the location and radius of the support of $m$.

\begin{figure}
	\centering
	\subfloat[]{\includegraphics[scale=0.2625]{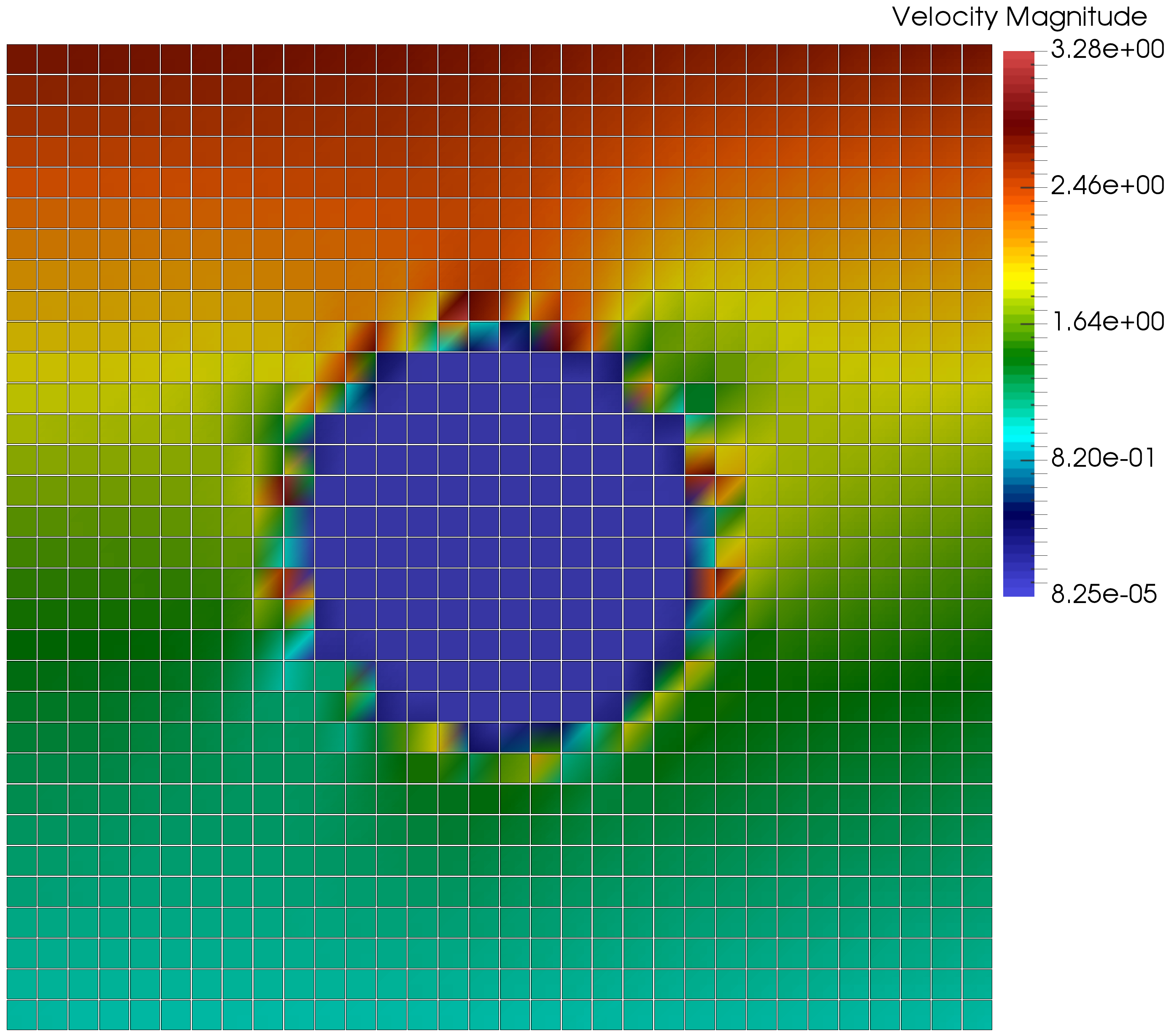}}
	\hfill
	\subfloat[]{\includegraphics[scale=0.2625]{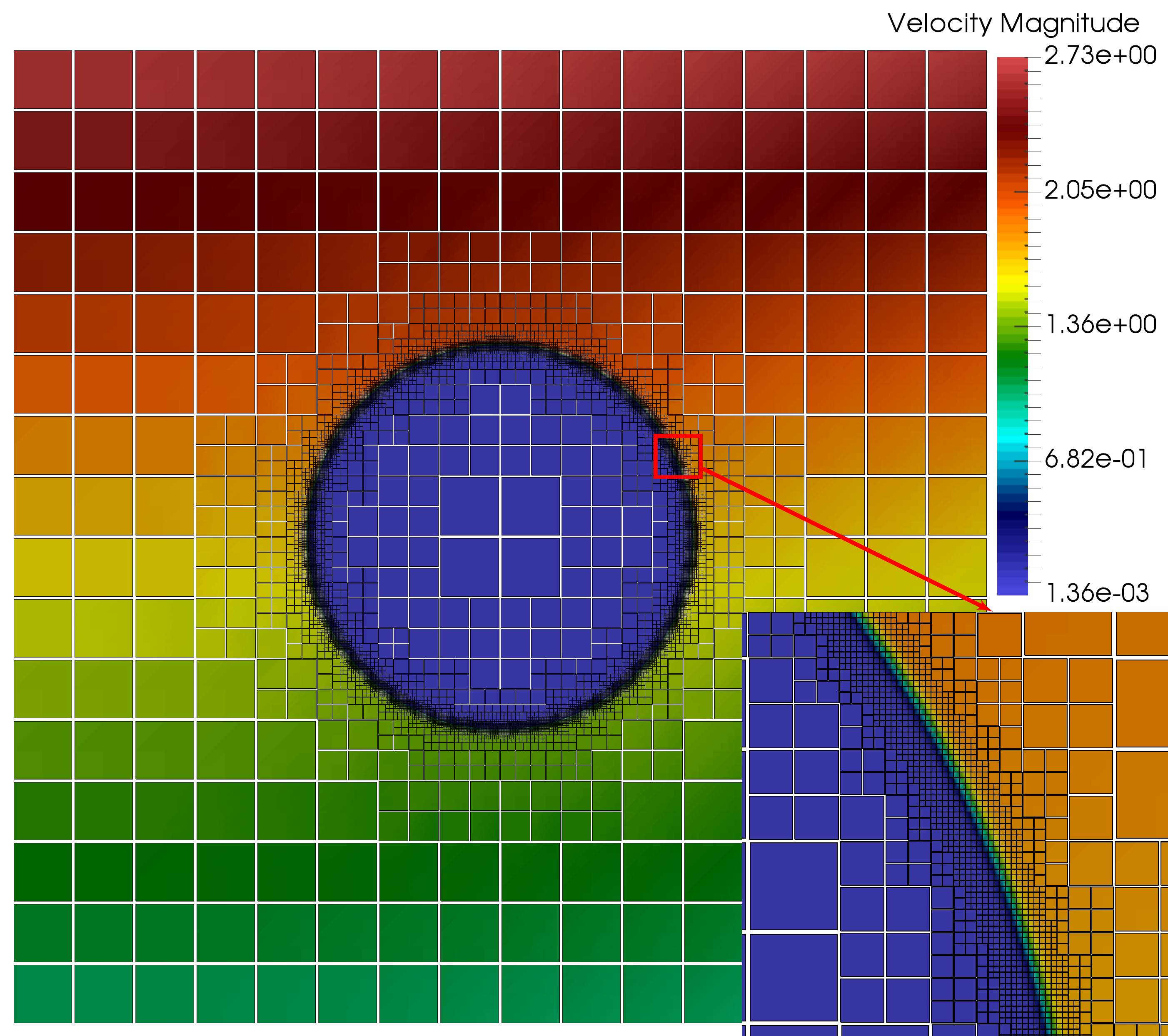}}
	\caption{Velocity magnitude of the numerical solution to
		the example defined in \secref{example4}, two dimensions, for a
		uniform level 5 mesh (a) and an adaptive mesh from level 4 to 10 (b).
                The artifacts arising with the uniform mesh are clearly visible.}
	\label{fig:err_e3_2d_vel}
\end{figure}

\begin{figure}
	\centering
	\includegraphics[scale=0.25]{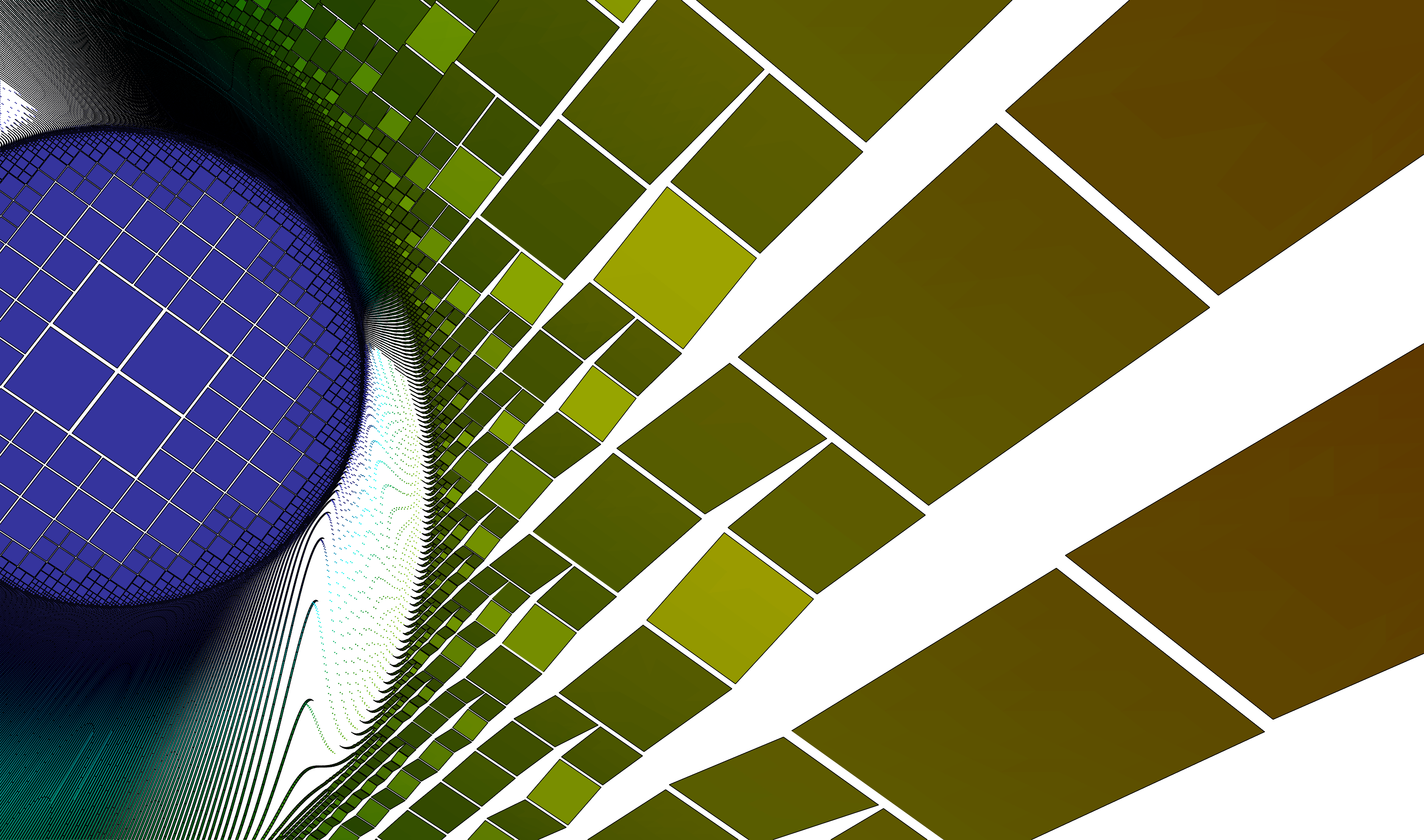}
	\caption{$y$-velocity magnitude extrusion illustrating
		a two dimensional \rtzero vector field for a level 4 to 10
		adaptively refined mesh. By construction of the Raviart-Thomas space,
                the $y$-velocity component is continuous in the $y$-direction
                and discontinuous in the $x$-direction.}
	\label{fig:err_e3_2d_wrap}
\end{figure}

The parameter $c$ allows us to tune the
coefficients to vary across the domain. From an application point of view, if $c$ is close
to one the function $m$ models a medium in which almost no flow is allowed within a
circle (sphere) centered at $\bm x_0$ with radius $a$. Given this information, it is clear
that the velocity field is likely to have a strong gradient within the ring-shaped
region $a < \norm{\bm{x}-\bm{x_0}}{} < b$.
We provide \figref{err_e3_2d_vel} and \figref{err_e3_2d_wrap} to illustrate this,
refining a given element on the mesh whenever it overlaps this region.
\begin{figure}
	\centering
	\subfloat[]{\includegraphics[scale=0.85]{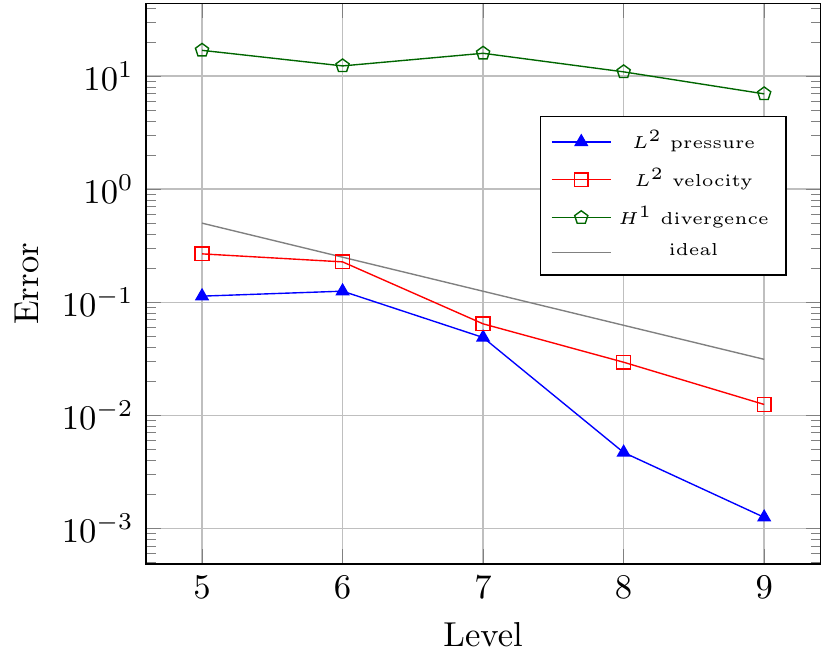}}
	\hfill
	\subfloat[]{\includegraphics[scale=0.85]{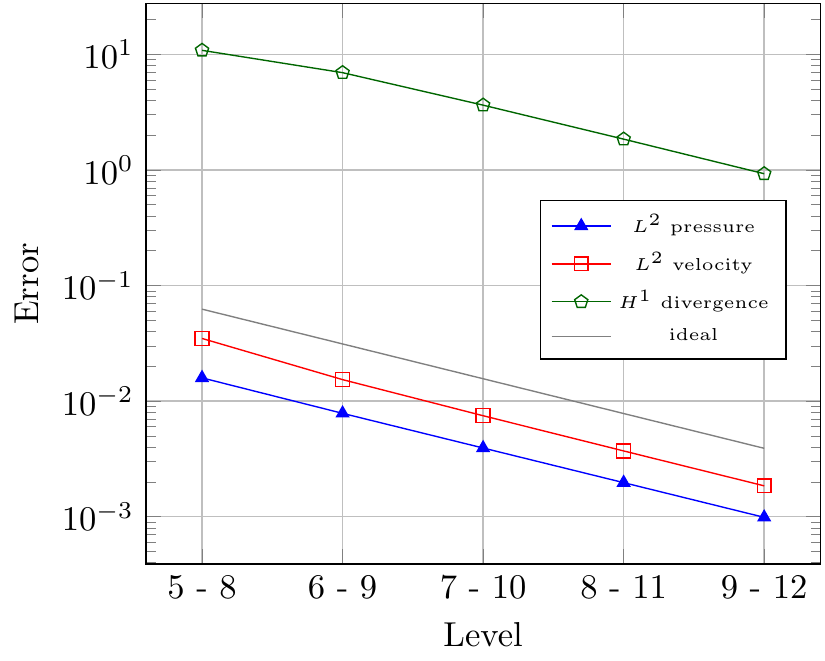}}
	\caption{Error plot for the numerical solution of a mixed Poisson system corresponding
	to the example defined in \secref{example4} (high conductivity contrast)
	in two dimensions for uniform (a) and adaptive (b) meshes.}
	\label{fig:err_e4_2d}
\end{figure}%
\begin{figure}
	\centering
	\subfloat[]{\includegraphics[scale=0.85]{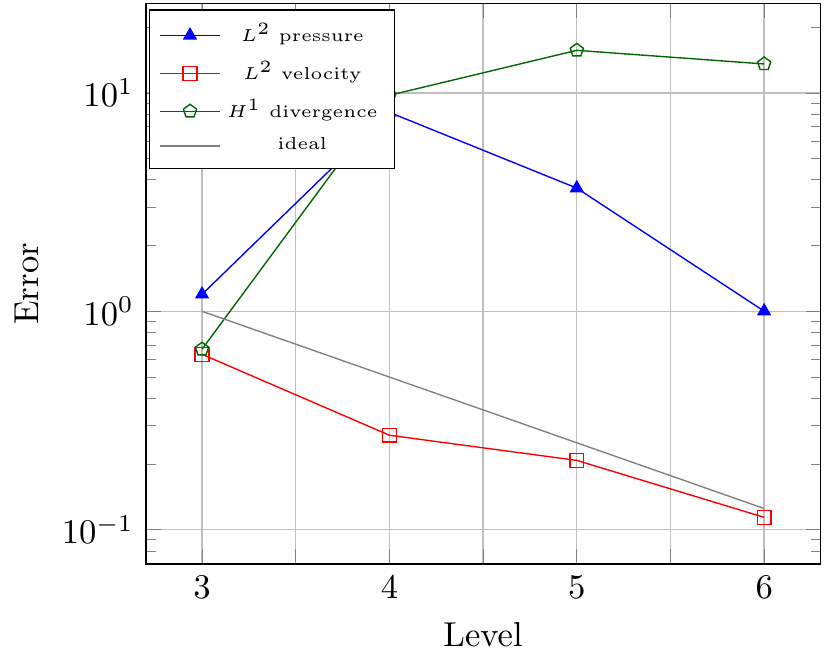}}
	\hfill
	\subfloat[]{\includegraphics[scale=0.85]{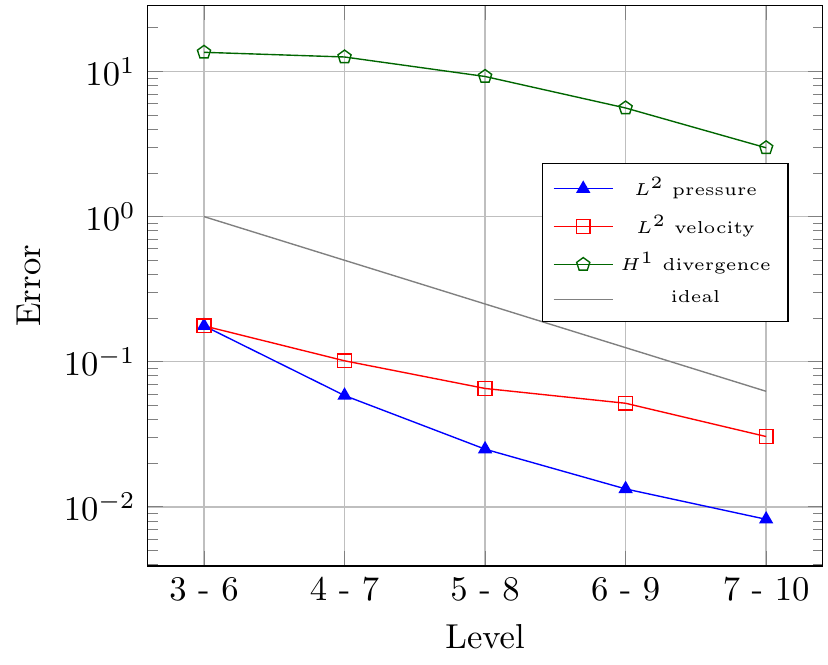}}
	\caption{Error plot for the numerical solution corresponding
		to the example defined in \secref{example4}
		(high contrast, $c=\mathrm{0.999}$) in three
		dimensions for uniform (a) and adaptive (b) meshes.}
	\label{fig:err_e4_3d}
\end{figure}%
% \begin{figure}
% 	\centering
% 	\includegraphics[scale=0.25]{e4_3d_vel_unif_zoom.pdf}
% 	\caption{Threshold plot from the computed velocity of a three dimensional
% 		version of the example defined in \secref{example4}.
% 		An adaptive mesh from level 4 to 8 was used in the calculation.}
% 	\label{err_e3_3d_vel}
% \end{figure}
\begin{figure}
	\centering
	\subfloat[]{\includegraphics[scale=0.85]{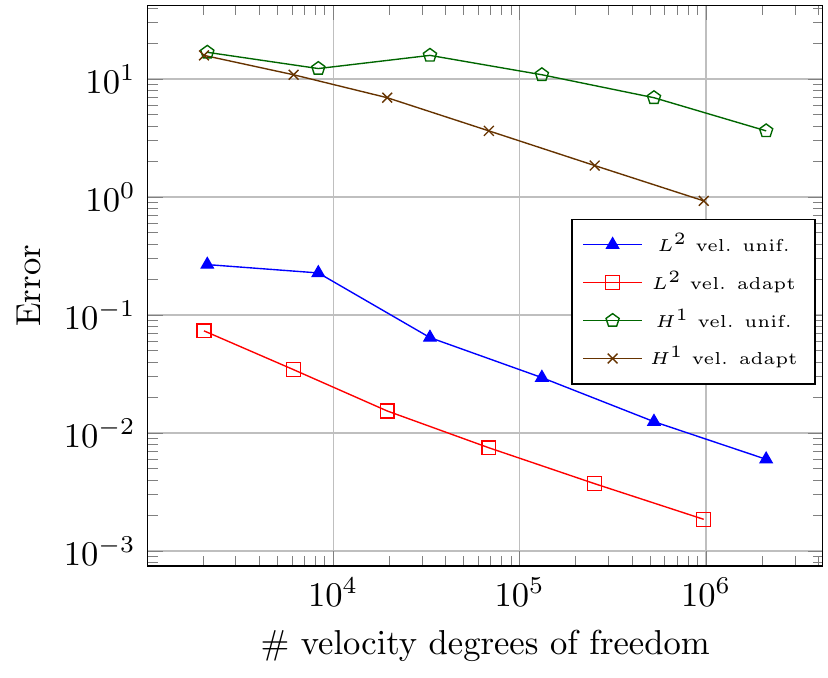}}
	\hfill
	\subfloat[]{\includegraphics[scale=0.85]{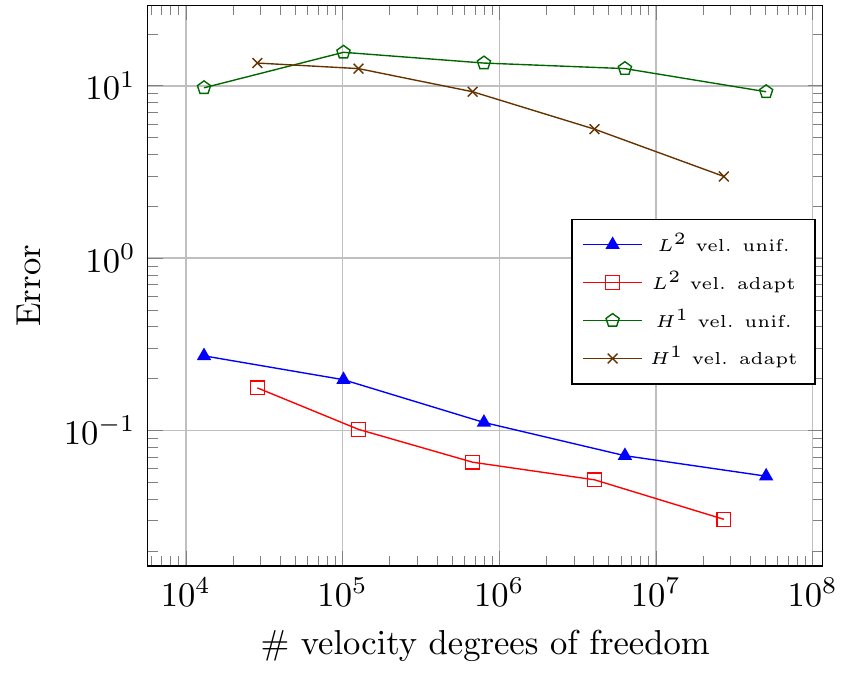}}
	\caption{For example defined in \secref{example4}, we compare the number of
		degrees of freedom with uniform and adaptive meshes against the
		$L^2$ and $H^{1}$ errors of the velocity for a two (a) and tree dimensional
		problem (b).}
	\figlab{e4_dof_vs_error}
\end{figure}%
\begin{table}
	\renewcommand{\tabcolsep}{1ex}
	\renewcommand{\arraystretch}{1.2}
	\centering
	\subfloat[]{\scalebox{0.65}{%
	\begin{tabular}{ ccccccc }
		\hline
		Level & \multicolumn{6}{c}{\# Iterations} \\ \cline{2-7}
		& NoPC 	& Diag 	& Schur & \multicolumn{3}{c}{SPAMG}  \\
		&	&	& 	& Uzawa	& Vanka one	& Vanka scale \\ \cline{5-7}
		$4$	&$>$1000&$>$1000& 36	& 12	& 8	& 9	\\
		$5$	& -	& -  	& 56	& 14	& 9	& 8 	\\
		$6$	& - 	& -	& 55	& 14	& 9     & 10	\\
		$7$	& - 	& -	& 47	& 15	& 10    & 11	\\
		$8$	& - 	& -	& 43	& 15	& 10    & 11	\\
		$9$	& - 	& -	& 41	& 15	& 11    & 11	\\
		\hline
	\end{tabular}
	}}
	\hfill
	\subfloat[]{\scalebox{0.65}{%
	\begin{tabular}{ ccccccc }
		\hline
		Level & \multicolumn{6}{c}{\# Iterations} \\ \cline{2-7}
		& NoPC 	& Diag 	& Schur & \multicolumn{3}{c}{SPAMG}  \\
		&	&	& 	& Uzawa	& Vanka one	& Vanka scale  \\ \cline{5-7}
		$4-7$   &$>$1000&$>$1000& 379 	& 12	& 8	& 8	\\
		$5-8$   & -	& -  	& 547	& 13	& 8	& 8 	\\
		$6-9$   & -	& - 	& 930	& 13	& 9     & 9	\\
		$7-10$  & -	& -	&$>$1000& 15	& 10    & 10	\\
		$8-11$  & -	& -	& - 	& 15	& 10    & 10	\\
		$9-12$	& -	& -	& - 	& 16	& 10    & 11	\\
		\hline
	\end{tabular}
	}}
	\caption{Number of iterations required by the GMRES solver for a two
		dimensional mixed Poisson system defined by the example in \secref{example4} discretized on a uniform (a) and adaptive mesh (b).}
	\tablab{ex4_iter_2d}
\end{table}%
\begin{table}
	\renewcommand{\tabcolsep}{1ex}
	\renewcommand{\arraystretch}{1.2}
	\centering
	\subfloat[]{\scalebox{0.85}{%
	\begin{tabular}{ ccccc }
		\hline
		Level & \multicolumn{4}{c}{\# Iterations} \\ \cline{2-5}
			& Schur & \multicolumn{3}{c}{SPAMG}  \\
			&	& Uzawa	& Vanka one  	& Vanka scale \\
		$3$   &	40	& 13	& 9	& 9	\\
		$4$   &	45	& 13	& 8	& 9	\\
		$5$   &	56 	& 15	& 9	& 9 	\\
		$6$   & 53  	& 15	& 10	& 10	\\
		$7$   & 47  	& 15	& 11	& 11 	\\
		\hline
	\end{tabular}
	}}
	\hfill
	\subfloat[]{\scalebox{0.85}{%
		\begin{tabular}{ ccccc }
			\hline
			Level & \multicolumn{4}{c}{\# Iterations} \\ \cline{2-5}
				& Schur & \multicolumn{3}{c}{SPAMG}  \\
				& 	& Uzawa	& Vanka one	& Vanka scale \\
			$3-6$   & 637  	& 15	& 11	& 11	\\
			$4-7$   & 855	& 14	& 10	& 11 	\\
			$5-8$   &$>$1000& 15	& 11	& 11	\\
			$6-9$   & -	& 16	& 12	& 12	\\
			$7-10$  & -     & 20	& 14	& 15	\\
			\hline
		\end{tabular}
	}}
	\caption{Number of iterations required by the GMRES solver for a three
		dimensional mixed Poisson system defined by the example in \secref{example4}
		discretized in a uniform (a) and adaptive mesh (b). We use the same setup as in
		\tabref{ex1_iter_2d}.}
	\tablab{ex4_iter_3d}
\end{table}%
In fact, adaptive refinement is necessary in this example to obtain optimal
convergence, which we show in \figref{err_e4_2d} (2D) and \figref{err_e4_3d} (3D).
The plot of accuracy versus degrees of freedom in \figref{e4_dof_vs_error}
supports this observation.

% \todo{Discuss how the Schur complement preconditioner does not work anymore.}

We display the iteration counts for various preconditioners in
\tabref{ex4_iter_2d} (2D) and \tabref{ex4_iter_3d} (3D).
While the Schur preconditioner still functions in 2D, it fails completely in
three space dimensions.
All three variants of SPAMG, on the other hand, lead to almost mesh
independent iteration counts.
The ``Vanka One'' variant seems best with an iteration count of just 10 to
reduce the relative error of the linear system of equations by 6 orders of
magnitude.
Thus, SPAMG is robust against coefficient functions whose magnitude varies by
at least a factor of 1000 in a narrow region.

\section{Conclusion}
\seclab{conclusion}

The purpose of this paper is to describe the SPAMG multigrid preconditioner for
saddle point systems.
It is a monolithic AMG method applicable to the block system arising from mixed
Poisson and Stokes discretizations.
One key element of the construction is a coupled prolongation operator that
stabilizes the Galerkin product from one level to the next coarser one.

Our numerical examples are obtained from a Raviart-Thomas discretization of the
mixed Poisson system.
We show the correctness of the solver and preconditioner by reproducing optimal
convergence rates and, more importantly, demonstrate mesh-independent iteration
counts of a preconditioned GMRES solve.
We see that the SPAMG preconditioner offers robustness against highly graded
adaptive meshes in both two and three space dimensions and tolerates strongly
varying coefficients of the PDE as well as a non-trivial, matrix-valued
conductivity.
In addition, we show that such is not achieved by a standard Schur complement
block preconditioner.

Essentially, mixed adaptive discretizations are handled well by SPAMG.
Some open questions remain, such as studying the robustness of the
preconditioner for even higher conductivity contrast on the order of $10^6$ and
extending it to largely disparate eigenvalues in the conductivity tensor (\ie
problems involving anisotropy).
Investigating the behavior of SPAMG for higher-order RT discretizations
is another possible extension.

%The numerical examples exposed in this document show that the SPAMG
%preconditioner yields iteration counts which are nearly independent
%of the mesh size. This behavior holds for uniform and locally refined meshes.
%This allows us to profit from the effectiveness of AMR: Less variables compared to
%the uniform mesh case lead to a better approximation of the expected solution
%(\eg \figref{e4_dof_vs_error}).
%Regarding the Schur complement strategy, our results show that in the presence of adaptive meshes denigrates the effectiveness of the preconditioner. Hence, SPAMG
%can be a valuable alternative as preconditioner for saddle point systems
%arising from a MFE discretization of a second order elliptic problems. In future work
%we would like to investigate the robustness of SPAMG with respect to more strongly heterogeneous
%coefficients and compare with newly developed strategies such as the auxiliary
%space multigrid presented in \cite{KrausLazarovLymberyEtAl16}.
%
%\todo{Try high order RT elements?}
%
%\todo{Use $\bm C=\bm D=\sqrt{\bm A}$ in BFBt ?}

\section*{Acknowledgments}

%\juelich

\trhcm

\sfb

\bibliographystyle{siam}
\bibliography{amg,group,ccgo_new}

\end{document}